\documentclass[12pt]{amsart}
\usepackage{amssymb}
\usepackage{amsmath}
\usepackage{longtable}
\newcommand\g{{\mathfrak g}}
\newcommand\h{{\mathfrak h}}
\newcommand{\p}{\mathfrak{p}}

\renewcommand{\a}{\mathfrak{a}}
\newcommand\m{\mathfrak m}
\renewcommand\k{\mathfrak k}
\newcommand\q{\mathfrak q}
\renewcommand\l{\mathfrak l}
\newcommand\z{\mathfrak z}

\renewcommand{\t}{\mathfrak{t}}
\newcommand\codim{\operatorname{codim}}
\newcommand\tr{\operatorname{tr}}
\newcommand\Spec{\operatorname{Spec}}
\newcommand\Quot{\operatorname{Quot}}

\newcommand\K{\mathbb K}
\newcommand\A{\mathbb A}

\newcommand\Z{\mathbb Z}

\newcommand\GL{\mathop{\rm GL}\nolimits}
\newcommand\Red{\operatorname{Red}}
\newcommand\SL{\operatorname{SL}}
\newcommand\Sp{\operatorname{Sp}}
\newcommand\Span{\operatorname{Span}}

\newcommand{\Ad}{\mathop{\rm Ad}\nolimits}

\newcommand{\rank}{\mathop{\rm rk}\nolimits}

\newcommand{\defe}{\mathop{\rm def}\nolimits}

\newcommand\Int{\mathop{\rm Int}\nolimits}

\newcommand{\im}{\operatorname{im}}
\newcommand{\Der}{\operatorname{Der}}
\newcommand\quo{/\!/}
\newtheorem{theorem}{Theorem}[subsection]
\newtheorem{proposition}[theorem]{Proposition}
\newtheorem{corollary}[theorem]{Corollary}
\newtheorem{lemma}[theorem]{Lemma}
\theoremstyle{definition}
\newtheorem*{specdefi}{Construction-definition of a compatible parabolic subgroup}
\newtheorem{example}[theorem]{Example}
\newtheorem{definition}[theorem]{Definition}
\newtheorem{remark}[theorem]{Remark}
  \numberwithin{equation}{section}
\numberwithin{theorem}{subsection} 
\title{Algebraic Hamiltonian actions
} 

\author{Ivan V. Losev}
\thanks{Supported by RFBR grant 05-01-00988}

\thanks
{Address: Department of Mathematics, Massachusetts Institute of Technology, 2-101
77 Massachusetts Avenue, Cambridge MA 02139, USA}
\thanks{e-mail: ivanlosev@math.mit.edu}
\begin{document}



\maketitle

\begin{abstract}
In this paper we deal with a Hamiltonian action of a reductive
algebraic group $G$ on an irreducible normal affine Poisson variety
$X$. We study the quotient morphism $\mu_{G,X}\quo G:X\quo
G\rightarrow \g\quo G$ of the moment map $\mu_{G,X}:X\rightarrow
\g$. We prove that  for a wide class of Hamiltonian actions
(including, for example, actions on generically symplectic
varieties) all fibers of the morphism $\mu_{G,X}\quo G$ have the
same dimension. We also study the "Stein factorization" of
$\mu_{G,X}\quo G$. Namely, let $C_{G,X}$ denote the spectrum of the
integral closure of $\mu_{G,X}^*(\K[\g]^G)$ in $\K(X)^G$. We
investigate the structure of the $\g\quo G$-scheme $C_{G,X}$. Our
results partially generalize those obtained by F. Knop for the
actions on cotangent bundles and symplectic vector spaces.
\end{abstract}

\section{Introduction}\label{SECTION_Intro}
\subsection{Main objects of study}\label{SUBSECTION_Intro1}
In this paper we study  Hamiltonian actions of reductive algebraic
groups on Poisson varieties. The ground field $\K$ is algebraically
closed and of characteristic 0.

A Poisson variety $X$ is a variety whose structure sheaf is a sheaf
of Poisson algebras. Poisson morphisms of Poisson varieties are
defined in an obvious way. Let $G$ be a reductive group acting on
$X$ by Poisson automorphisms. The action $G:X$ is said to be {\it
Hamiltonian} if it is equipped with a  $G$-equivariant linear map
$\g\rightarrow \K[X], \xi\mapsto H_\xi,$ where $\g$ denotes the Lie
algebra of $G$, such that the derivation $\{H_\xi,\cdot\}$ of the
algebra $\K(X)$ coincides with the velocity vector field $\xi_*$.
Under these conditions, the corresponding homomorphism
$S(\g)\rightarrow \K[X]$ is a homomorphism of Poisson algebras. A
Poisson variety equipped with a Hamiltonian action of $G$ is said to
be a Hamiltonian $G$-variety.

The morphism $\mu_{G,X}:X\rightarrow \g^*$ defined by
$\langle\mu_{G,X}(x),\xi\rangle=H_\xi(x), x\in X, \xi\in\g,$ is
called the {\it moment map}. Since $G$ is reductive, the algebra
$\g$ possesses a nondegenerate symmetric bilinear $G$-invariant form
$(\cdot,\cdot)$. We may assume additionally that this form is
nondegenerate on any Lie algebra of a reductive subgroup of $G$. Fix
such a form and identify $\g^*$ and $\g$. So we can consider
$\mu_{G,X}$ as a morphism $X\rightarrow \g$. We also consider a
morphism $\psi_{G,X}:X\rightarrow \g\quo G$, the composition of
$\mu_{G,X}$ and the quotient morphism $\pi_{G,\g}:\g\rightarrow
\g\quo G$.

It is interesting to study a kind of "Stein factorization" for the
morphism $\psi_{G,X}$. Let $X$ be a normal irreducible Hamiltonian
$G$-variety. Denote the integral closure of $\psi_{G,X}^*(\K[\g]^G)$
in $\K(X)^G$ by $A$. This is a finitely-generated subalgebra of
$\K[X]^G$. Put $C_{G,X}=\Spec(A)$. This is a normal irreducible
affine variety. There are two natural morphisms: the $G$-invariant
dominant morphism $\widetilde{\psi}_{G,X}:X\rightarrow C_{G,X}$  and
the finite morphism $\tau_{G,X}:C_{G,X}\rightarrow \g\quo G$. We
remark that, by the construction of $C_{G,X}$, $G$ permutes
transitively connected components of a general fiber of
$\widetilde{\psi}_{G,X}$. The variety $C_{G,X}$ and the morphisms
$\widetilde{\psi}_{G,X}:X\rightarrow C_{G,X}$ and
$\tau_{G,X}:C_{G,X}\rightarrow \g\quo G$ are the main objects of
study in this paper.

One of motivations for this study comes from the theory of
Hamiltonian actions of compact Lie groups on smooth symplectic
manifolds. Namely, let $K$ be a connected compact Lie group and $X$
a symplectic $K$-manifold.  As above, one can define the moment map
$\mu:X\rightarrow\k$. Choose a  Weyl chamber $C\subset \k$ and
consider the continuous map $\psi:X\rightarrow C$ mapping $x$ to
$K\mu(x)\cap C$. This is an analog of $\psi_{G,X}$ in this
situation. The map $\psi$ has the following properties
(see~\cite{Kirwan}):
\begin{enumerate}
\item The image of $\psi$ is a convex polytope.
\item All fibers of $\psi$ are connected.
\end{enumerate}

However,  general fibers of $\psi_{G,X}$ are, in general, not
connected even for connected $G$. The action of $G=\SL_2$ on
$\K^2\oplus \K^2$ provides an example. Therefore it seems that the
right analog of $\psi$ in the algebraic situation is the morphism
$\widetilde{\psi}_{G,X}$. We will see in the sequel that the image
of $\widetilde{\psi}_{G,X}$ possesses nice properties (at least for
sufficiently good, for example, generically symplectic, affine
varieties $X$). We also conjecture that for such varieties $X$ all
fibers of $\widetilde{\psi}_{G,X}\quo G:X\quo G\rightarrow C_{G,X}$
are irreducible\footnote{After this paper was submitted Friedrich Knop found a
counterexample to this conjecture}. This would  imply the connectedness property for
all fibers of $\widetilde{\psi}_{G,X}$.

The second motivation of our study comes from Invariant theory. It
turns out that the subalgebra $\K[C_{G,X}]\subset \K[X]^G$ is
closely related to the subalgebra of all functions lying in the
center of the Poisson algebra $\K(X)^G$. In particular, if $\K(X)^G$
is commutative, then $C_{G,X}$ is closely related to $X\quo G$.

The idea to study $C_{G,X},\widetilde{\psi}_{G,X},\tau_{G,X}$
belongs to F. Knop. In~\cite{Knop1} he showed that if $G$ is
connected, $X_0$ is a smooth irreducible $G$-variety and $X=T^*X_0$,
then $C_{G,X}$ is an affine space and the morphism $\psi_{G,X}$ is
equidimensional. He also described the morphism $\tau_{G,X}$.
Namely, there is a subspace $\a_{G,X_0}$ in a Cartan subalgebra $\t$
of $\g$ and a subgroup $W_{G,X_0}$ in the quotient
$N_G(\a_{G,X_0})/Z_G(\a_{G,X_0})$ such that $C_{G,X}\cong
\a_{G,X_0}/W_{G,X_0}$ and $\tau_{G,X}$ is the morphism
$\a_{G,X_0}/W_{G,X_0}\rightarrow \g\quo G$ induced by the
restriction of functions from $\g$ to $\a_{G,X_0}$. Since $C_{G,X}$
is an affine space, the group $W_{G,X_0}$ is generated by
reflections. The subalgebra $\K[C_{G,X}]$ coincides with the center
of the Poisson algebra $\K[X]^G$. Recently, F.Knop obtained the
analogous results for linear Hamiltonian actions, see~\cite{Knop3}.
Moreover, in the case of cotangent bundles all fibers of
$\widetilde{\psi}_{G,X}$ are connected, see \cite{Knop_connected}.

In the case when $X_0$ is a smooth irreducible $G$-variety the group
$W_{G,X_0}$ is an important birational invariant of $X$. In a
subsequent paper we will apply some results and constructions of the
present paper to the problem of the computation of $W_{G,X_0}$.

\subsection{Statement of the results}\label{SUBSECTION_Intro2}

We need some definitions.

\begin{definition}\label{definition:1.2.1}
Let $G$ be an arbitrary algebraic group. A $G$-variety $X$ is called
{\it $G$-irreducible} if $G$ acts transitively on the set of
irreducible components of $X$.
\end{definition}

$X$ is $G$-irreducible iff $\K(X)^G$ is a field.

Next, we define important numerical invariants of a Hamiltonian
$G$-variety. For a $G$-variety $X$ let $m_G(X)$ denote the maximal
dimension of a $G$-orbit on $X$.

\begin{definition}\label{definition:1.2.2}
Let $X$ be a $G$-irreducible Hamiltonian  $G$-variety.  The {\it
rank} of $X$ is, by definition, the number
$\rank_G(X)=m_G(\overline{\im\mu_{G,X}})$. The  difference
$m_G(X)-\rank_G(X)$ is called the {\it lower defect} of $X$ and is
denoted by $\underline{\defe}_G(X)$. The  {\it upper defect} is the
dimension of $\overline{\im\psi_{G,X}}$, it is denoted by
$\overline{\defe}_G(X)$.
\end{definition}

\begin{theorem}\label{theorem:1.2.3}
Suppose $X$ is an affine $G$-irreducible Hamiltonian $G$-\!\!
variety. Then the codimension of any fiber of $\mu_{G,X}\quo
G:X\quo G\rightarrow \g\quo G$ in $X\quo G$ is not less than
$\underline{\defe}_G(X)$.
\end{theorem}

\begin{definition}\label{definition:1.2.4}
A $G$-irreducible Hamiltonian $G$-variety $X$ is called {\it
equidefectinal} if $\underline{\defe}_G(X)=\overline{\defe}_G(X)$.
In this case we call $\underline{\defe}_G(X)=\overline{\defe}_G(X)$
the defect of $X$ and denote it by $\defe_G(X)$.
\end{definition}
We will see in Subsection~\ref{subsection_Ham4} that if $X$ is
generically symplectic (Definition~\ref{definition:2.2.4}) or
$m_G(X)=\dim G$, then $X$ is equidefectinal.

If $X$ is equidefectinal, then all fibers of $\mu_{G,X}\quo G$ have
the same dimension.  It is not clear whether $\psi_{G,X}$ possesses
this property. However, it is so when $X$ is smooth and symplectic.
We will prove this in a subsequent paper \cite{fibers}. Moreover, it can be shown
that any fiber of $\psi_{G,X}$ has a component of the "right"
dimension.

The next two theorems describe the morphism
$\tau_{G,X}:C_{G,X}\rightarrow \g\quo G$. To state them  we
introduce some factorization
$\tau_{G,X}=\tau^1_{G,X}\circ\tau^2_{G,X}$.

Let $X$ be a normal irreducible equidefectinal Hamiltonian
$G$-variety. For a point $x\in X$ in general position we put
$L=Z_{G^\circ}(\xi_s)$, where $\xi=\mu_{G,X}(x)$. Note that $L$ is
defined uniquely up to $G^\circ$-conjugacy. Put $\l^{pr}=\{ \xi\in
\l| \z_\g(\xi_s)\subset \l\}$. It follows from results of
Subsections \ref{subsection_Red1},\ref{subsection_Red2}, that
$\mu_{G,X}^{-1}(\l^{pr})$ is a normal $N_G(L)$-irreducible variety.
Choose a component $Y\subset \mu_{G,X}^{-1}(\l^{pr})$. Later on we
will see that the closure of the image of the projection
$p:\mu_{G,X}(Y)\rightarrow \l/[\l,\l]\cong \z(\l)$ in $\g$ is an
affine subspace in $\z(\l)$ of dimension $\defe_G(X)$
(Proposition~\ref{proposition:4.4.1}, Remark~\ref{remark:5.2.3}). We
denote this subspace by $\a_{G,X}^{(Y)}$. The group
$W_{G,X}^{(Y)}=N_G(L,Y)/L$ acts  on $\a_{G,X}^{(Y)}$ by affine
transformations. It turns out (see Subsection~\ref{subsection_Red2})
that $\tau_{G,X}=\tau^1_{G,X}\circ\tau^2_{G,X}$, where
$\tau^2_{G,X}:C_{G,X}\rightarrow \a_{G,X}^{(Y)}/W_{G,X}^{(Y)}$ is a
finite dominant morphism  and
$\tau^1_{G,X}:\a_{G,X}^{(Y)}/W_{G,X}^{(Y)}\rightarrow \g\quo G$ is
the finite morphism corresponding to the restriction of functions
from $\K[\g]^G$ to $\a_{G,X}^{(Y)}$.

In the case when $X_0$ is a quasi-affine algebraic variety,
$X=T^*X_0$ and $G$ is connected the pair
$(\a_{G,X}^{(Y)},W_{G,X}^{(Y)})$ is $G$-conjugate to the pair
$(\a_{G,X_0},W_{G,X_0})$ established by Knop. Our construction of
$\a_{G,X}^{(Y)},W_{G,X}^{(Y)}$ is inspired by Vinberg's variant of
the definition of $\a_{G,X_0},W_{G,X_0}$ (see~\cite{Vinberg1}). Note
that  Vinberg's construction is implicitly contained
in~\cite{Knop2}. The analogous construction of
$\a_{G,X_0},W_{G,X_0}$ for general $X_0$ was obtained
in~\cite{Timashev}.

\begin{theorem}\label{theorem:1.2.5}
Let $X$ be an equidefectinal normal irreducible affine Hamiltonian
$G$-variety. Then the morphism $\widetilde{\psi}_{G,X}\quo G:X\quo
G\rightarrow C_{G,X}$ is open. In particular,
$\im\widetilde{\psi}_{G,X}=\im(\widetilde{\psi}_{G,X}\quo G)$ is an
open subset of $C_{G,X}$. Moreover, there is a normal
$W_{G,X}^{(Y)}$-variety $Z$ such that $C_{G,X}\cong Z/W_{G,X}^{(Y)}$
and a finite morphism $\tau:Z\rightarrow \a_{G,X}^{(Y)}$ such that:
\begin{enumerate}
\item $\tau^2_{G,X}=\tau/W_{G,X}^{(Y)}$.
\item The morphism $\tau$ is \'{e}tale in all points of $\pi_{W_{G,X}^{(Y)},Z}^{-1}(\im
\widetilde{\psi}_{G,X})$.
\end{enumerate}
\end{theorem}

In particular, we get a partial description of singularities of
$\im\widetilde{\psi}_{G,X}$. For $Z$ we take the variety $C_{L,R}$,
where  $R$ is an affine normal irreducible equidefectinal
Hamiltonian $L$-variety constructed from $X$ (the {\it
$P_u$-reduction} of $X$, see below).

Under some additional restriction on the action  $G:X$ a more
precise statement can be obtained. The restriction is a presence of
some "good" action of the one-dimensional torus $\K^\times$ on $X$.
Let us give the precise definition.

\begin{definition}\label{definition:1.2.6}
An affine Hamiltonian $G$-variety $X$ equipped with an action
$\K^\times:X$ commuting with the action of $G$ is said to be {\it
conical} if the following two conditions are fulfilled
\begin{itemize}
\item[(Con1)] The  morphism $\K^\times\times X\quo G\rightarrow X\quo
G$ induced by the action $\K^\times: X$ can be extended to a
morphism $\K\times X\quo G\rightarrow X\quo G$.
\item[(Con2)] There exists a positive integer $k$ such that
$\mu_{G,X}(tx)=t^k\mu_{G,X}(x)$ for all $t\in \K^\times, x\in X$.
\end{itemize}
An integer $k$ satisfying the assumptions of (Con2) is called the
{\it degree} of $X$.
\end{definition}

For example, cotangent bundles and symplectic vector spaces with the
natural actions of $\K^\times$ are conical (see
Subsection~\ref{subsection_Ham3}).

\begin{theorem}\label{theorem:1.2.7}
Let $X$ be a conical Hamiltonian $G$-variety satisfying the
assumptions of Theorem~\ref{theorem:1.2.5}. Then $\a_{G,X}^{(Y)}$ is
a subspace in $\z(\l)$ and $\tau^1_{G,X}$ is an isomorphism. If, in
addition, $X$ is generically symplectic, then $\K[C_{G,X}]$
coincides with the subalgebra of all regular $G$-invariants lying in
the center of the Poisson field $\K(X)^G$.
\end{theorem}

Under the assumptions of the previous theorem, $\K[C_{G,X}]$
coincides with the center of $\K[X]^G$ provided
$\K(X)^G=\Quot(\K[X]^G)$. It turns out that the latter is true under
some additional assumptions.

\begin{definition}\label{definition:1.2.8} A Hamiltonian $G$-variety $X$ is called
 {\it strongly equidefectinal}
if there exists a stratification  $X=\coprod_i X_i$ by locally
closed $G$-irreducible equidefectinal Hamiltonian subvarieties (see
Definition~\ref{definition:3.1.4}) $X_i\subset X$ such that
$X_i=X_i^{max}$.
\end{definition}

Some classes of strongly equidefectinal Hamiltonian varieties are
listed in Subsection~\ref{subsection_Ham4}.


\begin{theorem}\label{theorem:1.2.9}
Suppose $X$ is a strongly equidefectinal normal affine irreducible
Hamiltonian $G$-variety. Then
\begin{enumerate}
\item A fiber of $\pi_{G,X}$ in general position contains a
dense $G$-orbit or, equivalently, $\K(X)^G=\Quot(\K[X]^G)$.
\item The following conditions are equivalent:
\begin{itemize}
\item[(a)] The action $G:X$ is stable, i.e. a fiber of $\pi_{G,X}$ in general position consists of one orbit.
\item[(b)] The stabilizer in general position for the action $G:X$ exists and is reductive.
\item[(c)] The subset of $\overline{\im\mu_{G,X}}$ consisting of
semisimple elements is dense in $\overline{\im\mu_{G,X}}$.
\end{itemize}
\end{enumerate}
\end{theorem}

\subsection{Some key ideas}\label{SUBSECTION_Intro3}
There are three main ingredients of the proofs. Let us give their
short (and not very precise) descriptions.

The first ingredient is the structure theory of a special class of
Hamiltonian $G$-varieties, namely central-nilpotent ones.

\begin{definition}\label{definition:1.3.1}
A Hamiltonian $G$-variety $X$ is called {\it central-nilpotent} (or,
shortly, CN) if  $\mu_{G,X}(x)_s\in \z(\g)$ for any $x\in X$.
\end{definition}

It is not very difficult to prove
Theorems~\ref{theorem:1.2.3},\ref{theorem:1.2.5} in the CN case.
Furthermore, irreducible normal affine CN Hamiltonian $G$-varieties
have a nice description provided $G$ is connected. Let us state this
result.

There are two important classes of affine CN Hamiltonian
$G$-varieties. Firstly, one can consider a Hamiltonian
$G/(G,G)$-variety $X_0$ as a Hamiltonian $G$-variety. Such
Hamiltonian $G$-varieties are clearly CN. To obtain one more
example, consider a nilpotent element $\eta\in\g$. By Example
\ref{example:3.2.7}, $X_1:=\Spec(\K[G/(G_\eta)^\circ])$ is a
Hamiltonian $G$-variety. This variety is again CN. Thus the product
$X_0\times X_1$ is CN. Consider a finite group $\Gamma$ acting on
$X_0\times X_1$ by Poisson automorphisms preserving the moment map.
We get a CN Hamiltonian variety $X_0\times X_1/\Gamma$. It turns out
that any affine irreducible normal CN Hamiltonian $G$-variety has
such a form provided $G\cong Z(G)^\circ\times (G,G)$. The last
requirement is not restrictive because any connected reductive
algebraic group possesses a covering satisfying this requirement.
Using this classification it is not very difficult to prove
Theorem~\ref{theorem:1.2.9}.

The second ingredient is the local theory of Hamiltonian actions on
quasi-projective varieties based on the Guillemin and Sternberg
local cross-section theorem (see~\cite{GS},\cite{Knop4}). Roughly
speaking,  the theorem reduces the study of a Hamiltonian
$G$-variety in an \'{e}tale neighborhood of a point $x\in X$ to the
study of a Hamiltonian action of the Levi subgroup
$Z_{G^\circ}(\mu_{G,X}(x)_s)$ on some locally closed subvariety of
$X$. This subvariety is called a {\it cross-section}. An example of a
cross-section is the Hamiltonian $L$-variety $Y\subset
\mu_{G,X}^{-1}(\l^{pr})$ mentioned above. Using local cross-sections
we complete  the proof of Theorem~\ref{theorem:1.2.9}.

To describe the third ingredient suppose that $X$ is an irreducible
normal affine equidefectinal Hamiltonian $G$-variety. Roughly
speaking, $Y$ is a CN Hamiltonian variety "approximating" $X$.
However, there is another CN Hamiltonian $L$-variety approximating
$X$ even better. This variety is constructed from $Y$ and an
appropriate parabolic subgroup $P\subset G$ with Levi subgroup $L$
and is called the {\it $P_u$-reduction of $X$ associated with $Y$}.
The name is chosen because our construction is, in some sense, a
modification of the Marsden-Weinstein reduction, see~\cite{MW}. The
idea is as follows. We want to consider the Marsden-Weinstein
reduction for the action $P_u:X$, that is, the quotient
$\mu_{G,X}^{-1}(\p)\quo P_u$. However, this quotient, in general,
seems to be very bad, possibly, it is not even a variety. Therefore
we take a "good" component of $\mu_{G,X}^{-1}(\p)$, namely
$Z=\overline{P_uY}$, and consider not the whole algebra
$\K[Z]^{P_u}$ but its subalgebra $A_Z$ generated by
$H_\xi|_Z,\xi\in\l,$ and $f|_Z$, $f\in \K[X]^{P_u}$. It turns out
that this subalgebra is finitely generated. The $P_u$-reduction $R$
is the normalization of the spectrum of the subalgebra. It is that
variety mentioned after Theorem~\ref{theorem:1.2.5}. $R$ possesses
the natural structure of a Hamiltonian $L$-variety (the hamiltonians
are $H_\xi|_{\overline{P_uY}}, \xi\in \l$). The $P_u$-reduction is
used to reduce the proofs of
Theorems~\ref{theorem:1.2.3},\ref{theorem:1.2.5} to the CN case.

\section{Poisson  varieties}
In Subsection~\ref{subsection_Poisson1} we define a Poisson (not
necessarily smooth) variety. In Subsection~\ref{subsection_Poisson2}
we define the Poisson bivector of a Poisson variety and study its
properties. In Subsection~\ref{subsection_Poisson3}  main examples
of Poisson varieties are given.  In
Subsection~\ref{subsection_Poisson4} we introduce a
stratification of a Poisson variety by smooth Poisson subvarieties
with the Poisson bivector of constant rank. Almost all definitions
and results of this section are well-known in the symplectic case.

\subsection{The main definition}\label{subsection_Poisson1} A commutative associative algebra $A$
with unit is called {\it Poisson} if it is equipped with a
skew-symmetric bilinear bracket $\{\cdot,\cdot\}:A\otimes
A\rightarrow A$ satisfying the Leibnitz and Jacobi identities, that
is
\begin{equation}\label{Leibniz}\{f,gh\}=\{f,g\}h+\{f,h\}g, \forall f,g,h\in
A,\end{equation}
\begin{equation}\label{Jacobi}
\{f,\{g,h\}\}+\{g,\{h,f\}\}+\{h,\{f,g\}\}=0, \forall f,g,h\in A.
\end{equation}
Thus the map  $f\mapsto \{f,g\}$ is a derivation of $A$ for any
$g\in A$.

Poisson homomorphisms of Poisson algebras are defined in a natural
way. An ideal $I\subset A$ is called Poisson if $\{A,I\}\subset I$.
For such an ideal $I$ the algebra $A/I$ possesses a unique Poisson
bracket such that the projection $A\rightarrow A/I$ is a Poisson
homomorphism.

\begin{proposition}\label{proposition:2.1.1}
Let $A\subset B$ be an algebraic extension of integral domains.
Suppose $A$ is equipped with a Poisson bracket $\{\cdot,\cdot\}_A$
and $B$ with some bracket $\{\cdot,\cdot\}_B$ (here a bracket is a
skew-symmetric bilinear operation satisfying the Leibnitz identity)
such that
\begin{equation}\label{eq:1.11}\{f,g\}_B=\{f,g\}_A, \forall f,g\in
A.
\end{equation} Then $\{\cdot,\cdot\}_B$ is a Poisson bracket. If brackets $\{\cdot,\cdot\}_B^1,\{\cdot,\cdot\}_B^2$ on $B$
satisfy (\ref{eq:1.11}), then they coincide.
\end{proposition}
\begin{proof}
This is a  consequence of the fact that $\{\cdot,\cdot\}$ is a
biderivation and the uniqueness of a lifting of a derivation for
algebraic extensions of integral domains (see~\cite{Leng}, Chapter
10). \end{proof}

\begin{definition}A variety $X$ is called {\it Poisson}, if its structure sheaf is a sheaf of Poisson algebras. A subvariety
$Y\subset X$ is called Poisson if its ideal sheaf  is a sheaf of
Poisson ideals. A morphism of Poisson varieties is called Poisson if
the corresponding homomorphisms of  algebras of sections of the
structure sheafs are Poisson.
\end{definition}

Note that a Poisson subvariety is naturally equipped with a
structure of a Poisson variety. Clearly, open subvariety of a
Poisson variety is Poisson.

Note that for any multiplicatively  closed subset $S$ of a Poisson
algebra $A$ the quotient algebra $A_S$ is equipped with a unique
Poisson bracket such that the natural homomorphism $A\rightarrow
A_S$ is Poisson. Thus a Poisson bracket on $\K[X]$ defines the
Poisson structure on $X$ provided $X$ is quasiaffine.

\begin{proposition}\label{proposition:2.1.2} Let $X$ be a Poisson variety. Then
\begin{enumerate}
\item Any irreducible component of $X$ is a Poisson subvariety.
\item Suppose $X$ is irreducible. The normalization
$\widetilde{X}$ of $X$ is equipped with a unique Poisson structure
such that the canonical morphism $\pi:\widetilde{X}\rightarrow X$ is
Poisson.
\end{enumerate}
\end{proposition}
\begin{proof}
We may assume that $X$ is affine. Recall that $\{f,\cdot\}$  is a
derivation of $\K[X]$ for any $f\in \K[X]$. Assertion 2 was proved by Kaledin
in \cite{Kaledin_normal}. Assertion 1
stems from the following lemma. \end{proof}

\begin{lemma}\label{lemma:2.1.4}
Suppose $A$ is a Noetherian $\K$-algebra. Minimal prime ideals  of
$A$ are stable under any derivation  $D\in \Der(A,A)$.
\end{lemma}
\begin{proof}
Localizing at a minimal prime ideal, we may assume that $A$ is a
local Artinian ring with the maximal ideal $\m$. Let $x\in\m$.
Choose an integer $n$ such that $x^{n-1}\neq x^n=0$. It remains to
note that $0=D(x^n)=nx^{n-1}Dx$. \end{proof}

\subsection{Poisson bivector}\label{subsection_Poisson2}
At first, we recall the relation between bivectors and 2-forms on
vector spaces. This material is standard, but we want to specify the
choice of signs.

Let $V$ be a finite dimensional vector space, $P\in \bigwedge^2V$.
The bivector $P$ induces the linear map $v:V^*\rightarrow V$ by
formula
\begin{equation}\label{eq:2.2:1}
\langle \alpha,v(\beta) \rangle=\langle P,\alpha\wedge\beta\rangle,
\alpha,\beta\in V^*.
\end{equation}

If  $P$ is nondegenerate, then  $v$ is an isomorphism. In general,
$P$ lies in $\bigwedge^2 v(V^*)$ and is a nondegenerate bivector in
this space. The map $v:V^*\rightarrow V$ is the composition of the
canonical surjection $V^*\rightarrow v(V^*)^*$, the isomorphism
$v(V^*)^*\rightarrow v(V^*)$ induced by $P\in \bigwedge^2 v(V^*)$
and the embedding $v(V^*)\hookrightarrow V$.

We can define the skew-symmetric nondegenerate bilinear form
$\omega_P$ on $v(V^*)$  by formula
\begin{equation}\label{eq:2.2:2}
\omega_P(v(\alpha),v(\beta))=\langle
P,\alpha\wedge\beta\rangle=\langle \alpha, v(\beta)\rangle
\end{equation}

Now let  $\omega$ be a nondegenerate skew-symmetric bilinear form on
$U\subset V$. We may consider $\omega$ as a nondegenerate bivector
in $\bigwedge^2 U^*$ and construct the bivector $P_\omega\in
\bigwedge^2U$ by formula (\ref{eq:2.2:2}). By embedding $\bigwedge^2 U$
into $\bigwedge^2 V$, we obtain the bivector in $\bigwedge^2 V$ with
$v(V^*)=U$. The maps $P\mapsto \omega_P$, $\omega\mapsto P_\omega$
are inverse to each other.

Now let $X$ be a variety. By a bracket on $X$ we mean a
skew-symmetric bilinear operation on $O_X$ satisfying the Leibnitz
identity. Brackets on $X$ are in one-to-one correspondence with
global bivectors, that is, global sections of the sheaf
$Hom_{O_X}(\bigwedge^2\Omega_X,O_X)$, where $\Omega_X$ is the sheaf
of K\"{a}hler differentials on $X$. If $X$ is smooth, then we get a
bivector in the usual sense. If a bracket satisfies (\ref{Jacobi}),
then the corresponding bivector is called {\it Poisson}. The bracket
corresponding to a bivector $P$ is given by
\begin{equation}\label{eq_bracket} \{f,g\}=P(df\wedge dg),
f,g\in \K(X).
\end{equation}

Now let $P$ be a bivector on $X$ and $x\in X$. Using the bivector
$P_x$,  we construct the linear map $v_x:T^*_xX\rightarrow T_xX$
defined by (\ref{eq:2.2:1}). Let $f$ be a rational function on $X$.
The vectors $v_x(df)$ form a vector field defined in the points of
the definition of $f$. This vector field is called the {\it
skew-gradient} of $f$, we denote it by $v(f)$. If $P$ is a Poisson
bivector, then, by the Jacobi identity for the bracket, the equality
$L_{v(f)}P=0$ holds, where $L$ denotes the Lie derivative.

Put \begin{equation}\label{eq:2.2:3}T^P_xX= \im v_x.\end{equation}
Clearly, $P_x\in \bigwedge^2 T^P_xX$. The bivector $P_x$ induces the
bilinear skew-symmetric nondegenerate form $\omega_x$ on $T^P_xX$ by
formula (\ref{eq:2.2:2}).

Let  $x\in X^{reg}$.   On the open subvariety $X^{max}\subset
X^{reg}$ consisting of all points $x$ such that $\rank
P_x=\max_{y\in X^{reg}}\rank P_y$ the spaces $T^P_xX$ form a locally
trivial vector bundle denoted by $T^PX$. We have the global section
$\omega$ of $\bigwedge^2 T^{P*}X$ over $X^{max}$ equal to $\omega_x$
in $x$. In the sequel we often call this section a 2-form. If
$X=X^{max}$, then $P$ is said to have constant rank.


\begin{definition}\label{definition:2.2.1}
Let $X$ be a smooth variety and $V$ a locally-trivial (in \'{e}tale
topology) subbundle of $TX$. The subbundle $V$ is called a {\it
distribution} on $X$. The distribution $V$ is called {\it
involutory}, if for any sections $\xi,\eta$ of $V$ on any \'{e}tale
neighborhood  of $X$ the commutator $[\xi,\eta]$ is also a section
of $V$.
\end{definition}

The  following proposition is standard (compare with \cite{CdSW},
Theorem 4.3, \cite{AG}, Subsection 3.2).

\begin{proposition}\label{proposition:2.2.2}
Let $X$ be a smooth variety and $P$ a bivector of constant rank on
$X$. Then the following conditions are equivalent
\begin{enumerate}
\item $P$ is Poisson.
\item The distribution $T^PX$ is involutory and  \begin{equation}\label{eq_sympl_form}
\begin{split}
&\omega([\xi,\eta],\zeta)+\omega([\eta,\zeta],\xi)+\omega([\zeta,\xi],\eta)=\\
&L_{\xi}(\omega(\eta,\zeta))+L_{\eta}(\omega(\zeta,\xi))+L_{\zeta}(\omega(\xi,\eta)).
\end{split}
\end{equation}
for all rational sections $\xi,\eta,\zeta$ of $T^PX$.
\end{enumerate}
Note that $P$ is uniquely determined by $T^PX$ and $\omega$.
\end{proposition}

\begin{definition}\label{definition:2.2.3}
A Poisson variety $X$ is called {\it symplectic}, if it is smooth
and for any $x$ from any irreducible component $X_0\subset X$ the
equality  $\rank_xP=\dim X_0$ holds.
\end{definition}

\begin{definition}\label{definition:2.2.4} An irreducible Poisson variety $X$ is said to be
{\it generically symplectic}, if $X^{max}$ is symplectic.
\end{definition}

If a Poisson variety $X$ is symplectic (resp.,  generically
symplectic), then $\omega$ is a symplectic form in the usual sense
on $X$ (resp., on $X^{max}$).

\subsection{Examples of Poisson varieties}\label{subsection_Poisson3}
\begin{example}\label{example:2.3.1}
Let $\g$ be an algebraic Lie algebra. The algebra $\K[\g^*]\cong
S(\g)$ possesses a unique Poisson bracket $\{\cdot,\cdot\}$ such
that $\{\xi,\eta\}=[\xi,\eta]$ for all $\xi,\eta\in \g$. Thus $\g^*$
is equipped with the structure of a Poisson variety. The
corresponding Poisson bivector $P$ is given by
$P_\alpha(\xi\wedge\eta)=\langle\alpha,[\xi,\eta]\rangle,
\alpha\in\g^*$. This implies $T^P_x=\g_*x$. A locally-closed
subvariety $X\subset\g^*$ is Poisson iff $X$ is $\Int(\g)$-stable.
In particular, an orbit  $O$ of the action $\Int(\g):\g^*$ is a
Poisson subvariety in $\g^*$. This variety is symplectic, the
corresponding symplectic form $\omega$ is called {\it the
Kostant-Kirillov form}. Explicitly,
$\omega_\alpha(\xi_*,\eta_*)=\alpha([\xi,\eta]), \alpha\in
O,\xi,\eta\in\g$.
\end{example}

\begin{example}\label{example:2.3.2}
Let $Y$ be a variety and ${\mathcal Vect}$ its sheaf of vector
fields. The vector bundle $T^*Y={\mathcal Spec}(S_{O_Y}({\mathcal
Vect}))$ is called the {\it cotangent bundle} of $Y$. The cotangent
bundle is locally trivial iff $Y$ is smooth. Let us equip $X=T^*Y$
with a natural Poisson structure. It is enough to do it locally and
check that the obtained structures are compatible. Thus one may
assume that $Y=\Spec(A)$ for an algebra $A$ of finite type. In this
case $\K[X]=S_A(D)$, where $D=\Der(A,A)$ is the module of
derivations of $A$. The algebra $S_A(D)$ possesses a unique Poisson
bracket $\{\cdot,\cdot\}$ such that
\begin{equation*}
\begin{split}
&\{f_1,f_2\}=0, \{f_1,d_1\}=d_1(f_1), \{d_1,d_2\}=[d_1,d_2]:=d_2d_1-d_1d_2,\\
& f_1,f_2\in A, d_1,d_2\in D.
\end{split}
\end{equation*}

The uniqueness follows from the fact that $A$ and $D$ generate
$S_A(D)$. Let us sketch the proof of the existence. Firstly, using
the construction of a tensor algebra, we prove that any $x\in D,a\in
A$ define derivations of the algebra $T_A(D)$ (commutators with
these elements). Then, considering $S_A(D)$ as a quotient of
$T_A(D)$, one can prove that $a,x$ define the derivations $d_a,d_x$
of $S_A(D)$. Using an analogous argument and the Leibnitz rule, we
can construct the derivation $d_f$ of $S_A(D)$ corresponding to
$f\in S_A(D)$. The bracket $\{g,f\}=d_f(g)$ has the required
properties.

The compatibility of these brackets follows from the uniqueness
property.

If $Y$ is smooth, the Poisson structure constructed above coincides
with the standard symplectic structure on $T^*Y$
(see.~\cite{Vinberg}, Ch. 2, Section 1.4).

Note that, by construction, any regular vector field on $Y$ defines
an element in $\K[T^*Y]$.
\end{example}

\begin{example}\label{example:2.3.3}
Let $X,Y$ be Poisson varieties. The product $X\times Y$ is naturally
equipped with a Poisson structure.
\end{example}

\begin{example}\label{example:2.3.4}
Let  $X$ be a Poisson variety, $Y$ a  variety, $\varphi:Y\rightarrow
X$ an \'{e}tale morphism. Let us show that $Y$ possesses a unique
Poisson structure such that $\varphi$ is a Poisson morphism. Let
$P_X$ be the Poisson bivector on $X$. There is a unique bivector
$P_Y$ on $Y$ such that $d\varphi P_Y=P_X$. It remains to show that
$P_Y$ is a Poisson bivector. The latter is an easy consequence of
Proposition~\ref{proposition:2.2.2}.
\end{example}

\begin{example}\label{example:2.3.5}
Let  $X$ be a Poisson variety, $Y$ a normal irreducible variety,
$\varphi:Y\rightarrow X$ a morphism. Suppose that $\varphi$ is \'{e}tale
in any point of an open subset $Y^0\subset Y^{reg}$ such that
$\codim_Y(Y\setminus Y^0)\geqslant 2$. Let us show that $Y$
possesses a unique Poisson structure such that $\varphi$ is a
Poisson morphism. By the previous example, $Y^0$ possesses a unique
Poisson structure such that $\varphi|_{Y^0}$ is a Poisson morphism.
Since $Y$ is normal, $\K[U]=\K[U\cap Y^0]$ for any open subset
$U\subset Y$. This allows one to define a Poisson structure on the
whole variety $Y$. Clearly, $\varphi$ is a Poisson morphism.
\end{example}

\subsection{Stratification of a Poisson variety}\label{subsection_Poisson4}
The following proposition appeared in \cite{Polishchuk}, Section 2. The proof is essentially contained
in Corollaries 2.3, 2.4.
\begin{proposition}\label{proposition:2.4.1}
Let $X$ be a Poisson variety. There exists a unique decomposition of
$X$ into the disjoint union of irreducible locally closed subvarieties
$X_i$ fulfilling the following conditions:
\begin{itemize}
\item[(a)] $\overline{X_i}$ is a Poisson subvariety of $X$.
\item[(b)] $X_i=\overline{X_i}^{max}$.
\end{itemize}
\end{proposition}
%
%
%


\section{Hamiltonian actions}
In Subsection~\ref{subsection_Ham1} we define Hamiltonian actions of
reductive groups on Poisson varieties and study their simplest
properties. In Subsections
\ref{subsection_Ham2},\ref{subsection_Ham3} we introduce some
examples of Hamiltonian (respectively, conical Hamiltonian)
varieties. In Subsection~\ref{subsection_Ham4} we describe some
classes of equidefectinal and strongly equidefectinal actions.
Finally, in Subsection~\ref{subsection_Ham5} we use Hamiltonian
actions of tori  to prove the Zariski-Nagata theorem on the purity
of branch locus.

In this section  $X$ is a Poisson variety (not necessarily smooth or
irreducible) and $G$ is a reductive group acting on  $X$ by Poisson
automorphisms.
\subsection{Main definitions and some
properties}\label{subsection_Ham1} Assume that there is a linear map
$\xi\mapsto H_\xi$  from $\g$ to  $ \K[X]$ satisfying the following
two conditions:
\begin{itemize}
\item[(H1)] $L_{\xi_*}f=\{H_\xi,f\}$ for any $f\in \K(X),\xi\in\g$.
\item[(H2)] The map $\xi\mapsto H_\xi$ is $G$-equivariant. 
\end{itemize}

\begin{definition}\label{definition:3.1.1}
An action $G:X$ together with a linear map $\xi\mapsto H_\xi$
satisfying (H1),(H2) is called {\it Hamiltonian}. If the action
$G:X$ is Hamiltonian, then $X$ is said to be a Hamiltonian
$G$-variety. The functions $H_\xi$ are called the {\it hamiltonians}
of the action. The morphism $\mu_{G,X}:X\rightarrow \g^*$ defined by
$\langle\mu_{G,X}(x),\xi\rangle=H_\xi(x)$ for all $x\in X,\xi\in\g$
is called the {\it moment map}.
\end{definition}

Since $\{H_\xi, H_\eta\}=L_{\xi_*}H_\eta=H_{[\xi,\eta]}$,
$\mu_{G,X}$ is a Poisson morphism.

When $X$ is symplectic, our definition coincides with the standard
one, see, for example,~\cite{Vinberg}, Ch.2, $\S2$.

In the sequel we fix a $G$-invariant  symmetric bilinear form
$(\xi,\eta)=\tr_V(\xi\eta)$ on $\g$, where $V$ is a locally
effective $G$-module. This form is nondegenerate on any subalgebra
$\h\subset \g$ corresponding to a reductive subgroup $H\subset G$.
We identify $\h^*$ with $\h$ using this form. So one may consider
the moment map as a morphism $X\rightarrow \g$.

Let  $\psi_{G,X}:X\rightarrow \g\quo G$ be the morphism defined by
$\psi_{G,X}=\pi_{G,\g}\circ\mu_{G,X}$.

\begin{remark}\label{remark:3.1.3} If  $H$ is a normal subgroup in $G$ and $X$ is a Hamiltonian $G/H$-variety,
then $X$ is naturally endowed with the structure of a Hamiltonian
$G$-variety. The moment map $\mu_{G,X}$ is the composition of
$\mu_{G/H,X}$ and the natural embedding $(\g/\h)^*\hookrightarrow
\g^*$. Conversely, if $X$ is a Hamiltonian  $G$-variety and a normal
subgroup $H\subset G$ acts trivially on $X$, we can consider $X$  as
a Hamiltonian $G/H$-variety with $\mu_{G/H,X}=\pi\circ\mu_{G,X}$,
where $\pi:\g\rightarrow \g/\h$ is a canonical projection. Note that
any fiber of $\psi_{G,X}$ is contained in a fiber of $\psi_{G/H,X}$.
\end{remark}

\begin{definition}\label{definition:3.1.4}
Let $X_1,X_2$ be Hamiltonian $G$-varieties, $\varphi:X_1\rightarrow
X_2$ a $G$-equivariant Poisson morphism. The morphism $\varphi$ is
said to be {\it Hamiltonian},
 if $\mu_{G,X_1}=\mu_{G,X_2}\circ \varphi$. If $Y$ is a $G$-stable subvariety
 of $X$ such that the embedding $Y\hookrightarrow X$ is
 Hamiltonian, then $Y$ is said to be a Hamiltonian subvariety of
 $X$.
\end{definition}


\begin{proposition}\label{proposition:3.1.7.5}
Let $X$ be a Hamiltonian $G$-variety and $x\in X^{max}$. For $v\in
T^P_xX,\xi\in\g,$ we have
\begin{equation*}\langle
d_x\mu_{G,X}(v),\xi\rangle=\omega_x(\xi_*x,v).\end{equation*} In
particular, $d_x\mu_{G,X}(T^P_xX)=\g_x^\perp$.
\end{proposition}
\begin{proof}
The proof is completely analogous to the symplectic case,
considered, for example, in \cite{Vinberg}, Chapter 2, Subsection
2.4. \end{proof}

We recall that the rank and the upper and  lower defects of a
$G$-irreducible Hamiltonian $G$-variety were defined in
Definition~\ref{definition:1.2.2}.

\begin{remark}\label{remark:3.1.8} If $X$ is symplectic, our
definition of the defect coincides with that given in
~\cite{Vinberg}, Subsection 2.5.
\end{remark}

The following properties of the rank and the defects  follow
directly from the definition.
\begin{lemma}\label{lemma:3.1.6}
\begin{enumerate}
\item Suppose that  $X$ is a $G$-irreducible Hamiltonian $G$-variety and $X_0$ is an irreducible component of $X$. Then
$\underline{\defe}_G(X)=\underline{\defe}_{G^\circ}(X_0)$,
$\overline{\defe}_G(X)=\overline{\defe}_{G^\circ}(X_0)$,
 $\rank_G(X)=\rank_{G^\circ}(X_0)$ .
\item Let $X_1,X_2$ be $G$-irreducible Hamiltonian $G$-varieties, $\varphi:X_1\rightarrow X_2$  a dominant generically finite
Hamiltonian morphism. Then
$\underline{\defe}_G(X_1)=\underline{\defe}_G(X_2)$,
$\overline{\defe}_G(X_1)=\overline{\defe}_G(X_2)$,
$\rank_G(X_1)=\rank_G(X_2)$.
\end{enumerate}
\end{lemma}

The next proposition is the main property of the upper and  lower
defects.
\begin{proposition}\label{proposition:3.1.7}
Let $X$ be  $G$-irreducible.  There are the inequalities
$$\underline{\defe}_G(X)\leqslant \overline{\defe}_G(X),
m_G(X)\leqslant \dim\overline{\mu_{G,X}(X)},$$ simultaneously
turning into equalities. If  $X$ is generically symplectic, then
these inequalities turn into equalities.
\end{proposition}
\begin{proof}
The case when $X$ is symplectic can be found in \cite{Vinberg},
Chapter 2, Subsections 2.4,2.6. In the general case the proof is
analogous. \end{proof}

\begin{remark}\label{remark:3.1.9} Let us give an example when the inequalities of Proposition \ref{proposition:3.1.7}
are strict. Let $X$ be an irreducible Poisson variety with the zero
Poisson bracket. Let the torus $(\K^\times)^n$ act trivially on $X$.
For the moment map one may take an arbitrary nonconstant morphism $X\rightarrow
\K^n$.
\end{remark}

\begin{corollary}\label{corollary:3.1.10}
The following assertions are equivalent: \begin{enumerate}
\item $\underline{\defe}_G(X)=\overline{\defe}_G(X)=\rank G$.
\item $m_G(X)=\dim G$.
\end{enumerate}
\end{corollary}
\begin{proof}
{$(1)\Rightarrow (2)$.} By Proposition~\ref{proposition:3.1.7},
$\dim \overline{\mu_{G,X}(X)}\quo G=\rank G$. Since a general fiber of
$\pi_{G,\g}$ is a single orbit, it follows that
$\overline{\mu_{G,X}(X)}=\g$. But
$\underline{\defe}_G(X)=\overline{\defe}_G(X)$.
Proposition~\ref{proposition:3.1.7} implies $m_G(X)=\dim G$.

The proof of $(2)\Rightarrow (1)$ is analogous. \end{proof}

\subsection{Examples of Hamiltonian varieties}\label{subsection_Ham2}
\begin{example}\label{example:3.2.1}
Let $X$ be a locally closed $G$-stable subvariety in $\g\cong\g^*$.
By Example~\ref{example:2.3.1},  $X$ is a Poisson variety. Put
$H_\xi=\xi|_X\in \K[X]$. It is checked directly that the pair
$X,(H_\xi)$ satisfies the conditions (H1),(H2). So $X$ is a
Hamiltonian  $G$-variety. The moment map is the embedding
$X\hookrightarrow\g$.
\end{example}

\begin{example}\label{example:3.2.2}
Let  $Y$ be a $G$-variety and $X=T^*Y$ (see
Example~\ref{example:2.3.2}). Being a global vector field on $Y$,
the velocity vector field $\xi_*,\xi\in \g,$ defines the function
$H_\xi\in \K[X]$. The group $G$ acts naturally on $X$. The pair
$X,(H_\xi)$ clearly satisfies (H2). Note that
$L_{\xi_*}\eta=[\xi_*,\eta]$ for any open subset $U\subset Y,
\eta\in {\mathcal Vect}(U),\xi\in\g$. It follows from the
construction of the Poisson structure that the pair $X, (H_\xi)$
satisfies (H1). The moment map is given by $\langle
\mu_{G,X}(y,\alpha),\xi \rangle=\langle \alpha,\xi_*y\rangle$, $y\in
Y,\alpha\in T^*_yY, \xi\in \g$.

If  $Y$ is smooth, the Hamiltonian structure constructed above
coincides with the standard one,  see~\cite{Vinberg}, Ch.2, $\S2$,
Example 1.
\end{example}

\begin{example}\label{example:3.2.3}
Let $X$ be a Hamiltonian $G$-variety, $\widetilde{X}$ its
normalization, $\varphi:\widetilde{X}\rightarrow X$ the canonical
morphism. The $G$-variety $\widetilde{X}$ can be equipped with a
unique Hamiltonian structure such that $\varphi$ is a Hamiltonian
morphism.
\end{example}

\begin{example}\label{example:3.2.4}
Let $X$ be a Hamiltonian $G$-variety, $Y$ a Poisson subvariety of
$X$. Since the ideal sheaf of $Y$ is stable under the brackets with
$H_\xi,\xi\in\g$,  the subvariety $Y\subset X$ is $G^\circ$-stable.
If $GY=Y$, then the pair $(Y,(H_\xi|_Y))$ satisfies (H1),(H2). Thus
$Y$ becomes a Hamiltonian subvariety of $X$.
\end{example}

\begin{example}\label{example:3.2.5}
Suppose $X_1,X_2$ are Poisson varieties, groups $G_1,G_2$ act on
$X_1,X_2$, respectively, and these actions are Hamiltonian. Then the
action $G_1\times G_2:X_1\times X_2$ is Hamiltonian. The moment map
is given by the formula $\mu_{G_1\times G_2, X_1\times
X_2}(x_1,x_2)=\mu_{G_1,X_1}(x_1)+\mu_{G_2,X_2}(x_2)$ for $x_1\in
X_1,x_2\in X_2$.
\end{example}

\begin{example}\label{example:3.2.6}
Let $X$ be a Hamiltonian $G$-variety, $Y$ an irreducible normal
$G$-variety and $\varphi:Y\rightarrow X$ a $G$-equivariant morphism
satisfying the assumptions of Example~\ref{example:2.3.5}. Then $Y$
is equipped with a unique Poisson structure such that $\varphi$ is a
Poisson  morphism. Since this structure is unique, it is
$G$-invariant. Let $H_\xi^X$ be the hamiltonians for the action
$G:X$. Put $H_\xi^Y=\varphi^*(H_\xi^X)$. Clearly, the pair
$Y,(H_\xi^Y)$ satisfies (H2). It is easily deduced from the
uniqueness of a lifting of a derivation (\cite{Leng}, Chapter 10)
that this pair satisfies also (H1). Note that
$\mu_{G,Y}=\mu_{G,X}\circ\varphi$ whence $\varphi$ is a Hamiltonian
morphism.
\end{example}

\begin{example}\label{example:3.2.7}
In particular, let  $\eta\in\g$,  $H$ be a subgroup of finite index
in $G_\eta$. Since all adjoint orbits have even dimension, the
algebra $\K[G\eta]$ is finitely generated (see.~\cite{VP}, Section
3.7). But $\K[G/H]$ is the integral closure of $\K[G\eta]$ in
$\K(G/H)$. Thus $\K[G/H]$ is finitely generated.  The natural
morphism $\varphi:\Spec(\K[G/H])\rightarrow \overline{G\eta}$
satisfies the assumptions of Example~\ref{example:3.2.6}. So
$X=\Spec(\K[G/H])$ is equipped with the structure of a Hamiltonian
$G$-variety. The equality $\mu_{G,X}=\varphi$ holds. Note that $G/H$
is an open subset in $\Spec(\K[G/H])$ and thus  a Hamiltonian
$G$-variety. This variety is symplectic.
\end{example}

\begin{example}\label{example:3.2.8}
There is an important special case of the previous construction. Let
$G=\Sp(V)$, where $V$ is a symplectic vector space with a constant
symplectic form $\omega$, and $\eta$  a highest weight vector of
$\g$. Then the $G$-variety $V$ coincides with
$\Spec(\K[G/(G_\eta)^\circ])$. The corresponding Poisson bivector
corresponds to $\omega$. The moment map is given by
(see~\cite{Vinberg}, Chapter 2, Example 2)
$\langle\mu_{G,V}(v),\xi\rangle=\frac{1}{2}\omega(\xi v,v), v\in
V,\xi\in\g$.
\end{example}

\begin{example}\label{example:3.2.9}
Let $X$ be a Hamiltonian $G$-variety with the hamiltonians $H_\xi,
\xi\in\g,$ and  $H$ a reductive subgroup in $G$. Then the
$H$-variety $X$ and the linear map $\h\rightarrow  \K[X],\xi\mapsto
H_\xi,$ satisfy  (H1),(H2).  Thus $H:X$ is a Hamiltonian action. The
moment map $\mu_{H,X}$ is the composition of $\mu_{G,X}$ and the
restriction map $\g^*\rightarrow \h^*$. In particular, any linear
action of a reductive group on a symplectic vector space becomes
Hamiltonian.
\end{example}

\begin{example}\label{example:3.2.10}
Suppose $X$ is a Hamiltonian $G$-variety, and a reductive group $H$
acts on  $X$ by Hamiltonian automorphisms. Suppose the good
categorical quotient $X\quo H$ exists (for example, $X$ is affine or
$X$ is quasiprojective and $H$ is finite).  Then $X\quo H$ is
equipped with a unique structure of a Hamiltonian $G$-variety such
that $\pi_{H,X}$ is a Hamiltonian morphism.
\end{example}

\subsection{Conical Hamiltonian varieties}\label{subsection_Ham3}
A conical Hamiltonian variety was defined in the Introduction,
Definition~\ref{definition:1.2.6}.

\begin{example}\label{example:3.3.1}
Let $H$ be a reductive  group. The group $\K^\times$ acts on $\h$ as
usual, that is $(t,x)\mapsto tx, t\in \K^\times, x\in \h$. Let $X$
be a closed  $H$-stable cone in $\h$ (a cone in $\h$, is, by
definition, a $\K^\times$-stable subset). For any reductive subgroup
$G\subset H$ the Hamiltonian $G$-variety $X$ (see
Examples~\ref{example:3.2.7} and~\ref{example:3.2.9}) equipped with
the action of $\K^\times$ induced from $\h$ is conical of degree 1.
\end{example}

\begin{example}\label{example:3.3.2}
Let $G:V$ be a linear Hamiltonian action (see
Examples~\ref{example:3.2.8}, \ref{example:3.2.9}). The Hamiltonian
$G$-variety $V$ together with the action $\K^\times:V$ given by
$(t,v)\mapsto tv$ is conical of degree 2.
\end{example}

\begin{example}\label{example:3.3.3}
Let $Y$ be an affine $G$-variety and $X=T^*Y$ (see
Example~\ref{example:3.2.2}). The variety $X$ is a vector bundle
over $Y$. Therefore there is the action $\K^\times:X$ by the
fiberwise multiplication. The Hamiltonian $G$-variety $X$ equipped
with this action of $\K^\times$ is conical of degree 1.
\end{example}

\begin{example}\label{example:3.3.4}
If $X$ is a conical Hamiltonian $G$-variety, then any $G$-stable
union of irreducible components of $X$ is a conical variety.
\end{example}

\begin{example}\label{example:3.3.5}
Let $X$ be a conical Hamiltonian $G$-variety. The action
$\K^\times:X$ can be lifted to the action of $\K^\times$ on the
normalization $\widetilde{X}$ of $X$. The Hamiltonian $G$-variety
$\widetilde{X}$ equipped with this action is conical.
\end{example}

The following lemma describes some basic properties of conical
varieties.
\begin{lemma}\label{lemma:3.3.7}
Let $X$ be a conical Hamiltonian $G$-variety of degree $k$. Then
\begin{enumerate}
\item $0\in \im\psi_{G,X}$.
\item Suppose $X$ is irreducible and normal. Then the subalgebra $\K[C_{G,X}]\subset
\K[X]^G$ (see the Introduction for the definition of $C_{G,X}$) is
$\K^\times$-stable. The morphisms
$\widetilde{\psi}_{G,X}:X\rightarrow C_{G,X},
\tau_{G,X}:C_{G,X}\rightarrow \g\quo G$ are $\K^\times$-equivariant,
where the action $\K^\times:\g\quo G$ is induced from the action
$\K^\times:\g$ given by $(t,x)\mapsto t^kx, t\in \K^\times, x\in\g$.
\item Under the assumptions of the previous assertion there exists a
unique point $\lambda_0\in C_{G,X}$ such that
$\tau_{G,X}(\lambda_0)=0$. For any point $\lambda\in C_{G,X}$ the
limit $\lim_{t\rightarrow 0}t\lambda$ exists and is equal to
$\lambda_0$.
\end{enumerate}
\end{lemma}
\begin{proof}
Let $x\in X$. It follows from (Con1) that there exists the limit
$y=\lim_{t\rightarrow 0}t\pi_{G,X}(x)$. By (Con2), $\mu_{G,X}\quo
G(y)=0$. This proves the first assertion.

 The morphism $\psi_{G,X}:X\rightarrow \g\quo G$ is
$\K^\times$-equivariant by (Con2).  Therefore the integral closure
of $\psi_{G,X}^*(\K[\g]^G)$ in $\K(X)^G$, that is, $\K[C_{G,X}]$, is
$\K^\times$-stable. This proves assertion 2.

Assertion 3 follows easily from the observations that $C_{G,X}$ is
irreducible and $\tau_{G,X}$ is $\K^\times$-equivariant.
 \end{proof}

\subsection{Equidefectinal and strongly equidefectinal Hamiltonian
varieties}\label{subsection_Ham4} The definitions of equidefectinal
and strongly equidefectinal Hamiltonian varieties were given in
the Introduction
(Definitions~\ref{definition:1.2.4},\ref{definition:1.2.8}).

\begin{lemma}\label{lemma:3.4.1}
Suppose $X$ is a $G$-irreducible Hamiltonian $G$-variety such that
 $m_G(X)=\dim G$ or some component of $X$ is generically symplectic. Then $X$ is equidefectinal.
\end{lemma}
\begin{proof}
If $m_G(X)=\dim G$ this follows directly from Proposition~\ref{proposition:3.1.7} and
Corollary~\ref{corollary:3.1.10}. Now suppose $X$ is generically symplectic.  By Proposition \ref{proposition:3.1.7}, $m_G(X)=\dim \overline{\im\mu_{G,X}}$. Note that
$\overline{\im\psi_{G,X}}=\overline{\im\mu_{G,X}}\quo G$. Every fiber of the quotient
morphism $\overline{\im\mu_{G,X}}\rightarrow \overline{\im\psi_{G,X}}$ consists
of finitely many orbits. Therefore the maximal dimension of a fiber coincides
with $\rank_G(X)$. Hence $m_G(X)-\rank_G(X)=\dim\overline{\im\psi_{G,X}}$.    \end{proof}

\begin{lemma}\label{lemma:3.4.2}
Suppose that $X$ is a Hamiltonian $G$-variety such that any stratum
described in Proposition~\ref{proposition:2.4.1} is a symplectic
variety. Then $X$ is strongly equidefectinal. In particular, any
symplectic variety $X$ is strongly equidefectinal.
\end{lemma}
\begin{proof}
This follows easily from Lemma~\ref{lemma:3.4.1}. \end{proof}

\begin{lemma}\label{lemma:3.4.3}
Let $X,Y$ be  Poisson varieties and $\varphi:X\rightarrow Y$ a
Poisson morphism. Then
\begin{enumerate}
\item For any Poisson subvariety $X'\subset X$ the closure of
$\varphi(X')$ is a Poisson subvariety of $Y$.
\item Suppose that the stratification of $Y$  consists of symplectic
varieties. Then the same is true for $X$ provided $\varphi$ is
finite and dominant.
\end{enumerate}
\end{lemma}
\begin{proof}
The first assertion is straightforward.

Proceed to  assertion 2. Note that any irreducible Poisson
subvariety of $Y$ is generically symplectic. Let $X'$ be an
irreducible locally closed Poisson subvariety of $X$. The subvariety
$\overline{\varphi(X')}\subset Y$ is Poisson. Since $\varphi$ is
finite, the ranks of Poisson bivectors on $X_1$ and on
$\overline{\varphi(X_1)}$ coincide. Thus $X'$ is generically
symplectic. It remains to apply this observation to the strata of
$X$. \end{proof}

\begin{corollary}
Let $G$ be a connected reductive group, $G_0$ its reductive
subgroup, $\eta\in\g$ and $H$ a subgroup of finite index in
$G_\eta$. The Hamiltonian $G_0$-variety $\Spec(\K[G/H])$ is strongly
equidefectinal.
\end{corollary}
\begin{proof}
Apply Lemma~\ref{lemma:3.4.3} to the morphism
$\Spec(\K[G/H])\rightarrow \overline{G\eta}$ and use
Lemma~\ref{lemma:3.4.2}. \end{proof}

The following lemma is used in the proof of
Theorem~\ref{theorem:1.2.9}.

\begin{lemma}\label{lemma:16.6}
Let $X,Y$ be Hamiltonian  $G$-varieties and $\varphi:X\rightarrow Y$
an \'{e}tale Hamiltonian morphism. If $Y$ is strongly equidefectinal,
then so is $X$.
\end{lemma}
\begin{proof}
Let $Y=\coprod_i Y_i$ be the stratification of $Y$ satisfying the
claims of Definition~\ref{definition:1.2.8}. Since $\varphi$ is
\'{e}tale, one can see that  $X_i=\varphi^{-1}(Y_i)$ is a Poisson subvariety
of $X$.  The morphism $\varphi|_{X_i}:X_i\rightarrow Y_i$ is
automatically Poisson and \'{e}tale. Thence $X_i$ is smooth. Let
$X_i=\coprod X_i^j$ be a unique stratification of $X_i$ by
$G$-irreducible unions of components. This stratification satisfies
the requirements of Definition~\ref{definition:1.2.8}.
\end{proof}

\subsection{An application of Hamiltonian actions: the Zariski-Nagata theorem on the purity of branch
locus}\label{subsection_Ham5} In this subsection we generalize the
Zariski-Nagata theorem, see, for example,~\cite{Danilov}, Chapter 4,
Subsection 1.4. However,  our generalization can be easily deduced
from this theorem.
\begin{proposition}\label{proposition:3.3}
Let $X,Y$ be irreducible varieties of the same dimension,
$\varphi:Y\rightarrow X$ a dominant morphism. Suppose
 $Y$ is normal, $X$ is smooth.  Then the complement in
 $Y$ to the open subset $Y^0=\{y\in Y| \varphi\text{ is \'{e}tale in } y\}$
 is a subvariety of pure codimension 1.
\end{proposition}
\begin{proof}
 Removing all components of $Y\setminus Y^0$ of codimension 1, we
may assume that $\codim_YY\setminus Y^0\geqslant 2$. Further, any
point $x\in X$ has a neighborhood admitting an \'{e}tale morphism to
$\A^n$ (see, for example,~\cite{Danilov}, Chapter 2, Subsection
6.3.). So we may assume $X=\A^n, Y$ is affine. Put
$T=(\K^\times)^n$. The action $T:T^*T\cong T\times \A^n$ is
Hamiltonian (Example~\ref{example:3.2.2}). Consider the action
$T:T\times Y$, where $T$ acts on $Y$ trivially and on $T$ by left
translations, and the $T$-equivariant morphism $\Phi:T\times
Y\rightarrow T\times \A^n, \Phi(t,y)=(t,\varphi(y))$. This morphism
satisfies the conditions of Example~\ref{example:3.2.6}, thus the
action $T:T\times Y$ is Hamiltonian. If $Y$ is smooth, then
$\varphi$ is \'{e}tale because $Y\setminus Y^0$ is the set of zeroes of
the Jacobian of $\varphi$. Assume that $Y$ is not smooth. Corollary 2.4 in \cite{Polishchuk} implies that
$T\times Y^{sing}=(T\times Y)^{sing}$ is a Poisson subvariety in
$T\times Y$. The action $T:(T\times Y)^{sing}=T\times Y^{sing}$ is
Hamiltonian (Example~\ref{example:3.2.4}). Choose an irreducible
component $Z\subset T\times Y^{sing}$. We see that
$\overline{\defe}_{T}(Z)\leqslant \dim Z\quo T=\dim Y^{sing}<\dim
Y=\dim T=\underline{\defe}_T(Z)$. This contradicts
Proposition~\ref{proposition:3.1.7}. \end{proof}

\section{Central-nilpotent Hamiltonian varieties}
Throughout this section $G$ is a connected reductive algebraic
group, $X$ is an  irreducible quasiprojective CN Hamiltonian
$G$-variety (see Definition~\ref{definition:1.3.1}).

Since $X$ is CN as a Hamiltonian $G$-variety, $X$ is CN also as a
$(G,G)$-variety.  In other words, the image of $\mu_{(G,G),X}$
consists of nilpotent elements. The number of nilpotent orbits is
finite.  Thus there is a unique open orbit $G\eta\subset
\overline{\im\mu_{(G,G),X}}$. Put $X^0=\mu_{G,X}^{-1}(G\eta)$. This
is a $G$-stable open subset of $X$.

In the first two subsections we study the structure of the
Hamiltonian $G$-variety $X$. In Subsection~\ref{subsection_CN1} we
describe the variety $X^0$ for an arbitrary normal variety $X$. In
Subsection~\ref{subsection_CN2} we describe the whole variety $X$
provided that $X$ is, in addition, affine. The description is
carried out for groups $G$ such that $Z(G)^\circ\cap (G,G)=\{1\}$,
in  this case $G\cong Z(G)^\circ\times (G,G)$. This requirement is
not restrictive:  any connected reductive group possesses a covering
satisfying the requirement.

In Subsection~\ref{subsection_CN3} we prove that the dimension of a
fiber of $\psi_{G,X}$ for an irreducible  CN Hamiltonian variety $X$
does not exceed $\dim X-\underline{\defe}_G(X)$. This gives us a
proof of Theorem~\ref{theorem:1.2.3} for CN varieties.

In Subsections~\ref{subsection_CN4} and \ref{subsection_CN5}
Theorems~\ref{theorem:1.2.5},\ref{theorem:1.2.9} are proved for
central-nilpotent varieties $X$.

\subsection{The structure of $X^0$}\label{subsection_CN1}
Let $X$ be normal, $G\cong Z(G)^\circ\times (G,G)$, $\eta,X^0$ such
as above. We put $O:=G/(G_\eta)^\circ\cong
(G,G)/((G,G)_\eta)^\circ.$ This is a Hamiltonian $G$-variety (see
Example~\ref{example:3.2.7}).

Let $X_0$ be a quasiprojective Hamiltonian $Z(G)^\circ$-variety.
Since $G=Z(G)^\circ$ $\times (G,G)$, we can consider $X_0$ as a
Hamiltonian $G$-variety. Let $\Gamma$ be a finite group acting on
$X_0\times O$ by Hamiltonian automorphisms. The quotient $(X_0\times
O)/\Gamma$ is equipped with  the natural Hamiltonian structure (see
Example~\ref{example:3.2.10}).

The main  result of this subsection is the following

\begin{theorem}\label{theorem:4.1.1} There exist a Hamiltonian
$Z(G)^\circ$-variety $X_0$ and a finite group $\Gamma$ acting freely
on $X_0\times O$ by Hamiltonian automorphisms such that $X^0\cong
(X_0\times O)/\Gamma$ (the isomorphism of Hamiltonian
$G$-varieties).
\end{theorem}

The theorem is proved in the following lemmas.

\begin{lemma}\label{lemma:4.1.2}
$((G,G)_{\mu_{G,X}(x)})^\circ\subset (G,G)_x\subset
(G,G)_{\mu_{G,X}(x)}$ for any $x\in X^0$.
\end{lemma}
\begin{proof}
By the choice of $\eta$,
$\overline{\im\mu_{(G,G),X}}=\overline{G\eta}$.
Proposition~\ref{proposition:3.1.7} implies that $\dim (G,G)x$
$=\dim (G,G)\eta$. Since $(G,G)_x\subset (G,G)_{\mu_{G,X}(x)}$, we
are done. \end{proof}

\begin{lemma}\label{lemma:4.1.3}
Let $G$ be an algebraic group, $H$ its subgroup, $Y$ a
quasiprojective $H$-scheme, $X=G*_HY$ (this homogeneous bundle
exists  by~\cite{VP}, Section 4.8). Then
\begin{enumerate}
\item $Y$ is $H$-irreducible iff $X$ is $G$-irreducible.
\item $X$ is normal iff $Y$ is so.
\end{enumerate}
\end{lemma}
\begin{proof}
The first assertion follows directly from the definition of
homogeneous bundles.

Proceed to assertion 2. Clearly, if $Y$ is a normal variety, then so
is $G*_HY$. Suppose now that $G*_HY$ is a normal variety.  Let
$Y_{red}$ denote the variety associated with $Y$. The canonical
morphism $G*_HY\rightarrow G*_HY_{red}$ is an isomorphism. Hence
$Y=Y_{red}$. If $\widetilde{Y}$ is the normalization of $Y$, then
$G*_N\widetilde{Y}$ is the normalization of $G*_NY$. Thus $Y$ is
normal. \end{proof}

\begin{lemma}\label{lemma:4.1.4}
The morphism of schemes
$\varphi:G*_{G_\eta}\mu_{(G,G),X}^{-1}(\eta)\rightarrow X$ is an
open embedding with the image $X^0$. The subscheme
$\mu_{(G,G),X}^{-1}(\eta)\subset X$ is a normal $G_\eta$-irreducible
subvariety.
\end{lemma}
\begin{proof}
Obviously $\im\varphi=X^0$. On the other hand,  the morphism
$\mu_{(G,G),X^0}:X^0\rightarrow G/G_\eta$ is $G$-equivariant. So
$\varphi$ induces an isomorphism of $G$-schemes $X^0\cong
G*_{G_\eta}\mu_{(G,G),X}^{-1}(\eta)$. By Lemma~\ref{lemma:4.1.3},
$\mu_{(G,G),X}^{-1}(\eta)$ is normal and $G_\eta$-irreducible.
\end{proof}

Choose a component $X_0$ of $\mu_{(G,G),X}^{-1}(\eta)$ and denote by
$H$ its stabilizer in the group $G_\eta$. Note that
$(G_\eta)^\circ\subset H$ and that $X^0\cong G*_{H}X_0$. The action
of $((G,G)_\eta)^\circ$ on $X_0$ is trivial by
Lemma~\ref{lemma:4.1.2}.

Put $\Gamma=H/(G_\eta)^\circ\cong (H\cap (G,G))/((G,G)_\eta)^\circ$.
This is a finite group acting freely on $O$ by Hamiltonian
automorphisms (the action is by the right translations). Since the
action of $((G,G)_\eta)^\circ$ on $X_0$ is trivial, $\Gamma$ acts
also on $X_0$ by $Z(G)^\circ$-automorphisms. Further, $\Gamma\cong
H/H^\circ$ acts on $G*_{H^\circ}X_0$,
$\gamma[g,x]=[g\widetilde{\gamma}^{-1},\widetilde{\gamma}x], g\in
G,x\in X_0,\gamma\in \Gamma$, where $\widetilde{\gamma}$ is an
element from $H$ mapping to $\gamma$ under the natural projection
$H\rightarrow H/H^\circ$. So the natural morphism
$G*_{H^\circ}X_0\rightarrow G*_{H}X_0\cong X^0$ is the quotient for
the action of $\Gamma$.

Consider the natural morphism $(G,G)\times X_0\rightarrow
G*_{H^\circ}X_0, (g,x)\mapsto [g,x]$. Since $((G,G)_\eta)^\circ$
acts trivially on $X_0$, this morphism factors through $O\times
X_0\cong (G,G)/((G,G)_\eta)^\circ\times X_0\rightarrow
G*_{H^\circ}X$. The latter is clearly an isomorphism. The action of
$\Gamma$ on $O\times X_0$ as on the product of $\Gamma$-varieties
coincides with the action on $G*_{H^\circ}X_0$. This action is free.

To complete the proof of Theorem~\ref{theorem:4.1.1} it remains to
prove the following

\begin{lemma}\label{lemma:4.1.5}
 There exists a
Poisson bracket on $X_0$ such that the action $Z(G)^\circ:X_0$ is
Hamiltonian with the moment map $\mu_{G,X}|_{X_0}-\eta$, the group
$\Gamma$ acts on $X_0$ by Hamiltonian automorphisms and the morphism
$O\times X_0\rightarrow X^0$ is Hamiltonian.
\end{lemma}
\begin{proof}
The morphism $\pi_{\Gamma,O\times X_0 }:O\times X_0\rightarrow  X^0$
is \'{e}tale. Lift the Hamiltonian structure from $X^0$ to
$G*_{H^\circ}X_0\cong O\times X_0$. Clearly, $\Gamma$ acts on
$O\times X_0$  by Hamiltonian automorphisms.

Let us introduce a Poisson bracket on $X_0$. For $x\in X_0$ the map
$\mu_{(G,G),X}$ is a covering $(G,G)x\rightarrow G\eta$
(Lemma~\ref{lemma:4.1.2}) thus $(G,G)x$ is equipped with the Poisson
bracket lifted from $G\eta$. We denote this bracket by
$\{\cdot,\cdot\}^{(G,G)x}$.

Identify $X_0$ with  $X_0\times \{eH^\circ\}\subset X_0\times O$.
For $f,g\in O_{O\times X_0,x}$ and $y\in X_0$ put
\begin{equation}\label{eq_bracket_X0}
\{f,g\}^{X_0}(y)=\{f,g\}(y)-\{f|_{(G,G)y},g|_{(G,G)y}\}^{(G,G)y}.
\end{equation}

$\{f,g\}^{X_0}$ is an element of $O_{X_0,x}$.  For $t\in Z(G)^\circ,
y\in X_0$ the equality $\{tf,tg\}^{X_0}(ty)=\{f,g\}^{X_0}(y)$ holds.
Let us check that $\{\cdot,\cdot\}^{X_0}$ is a Poisson bracket on
$\K(X_0)$. In the proof we may assume that $G$ is semisimple.

Denote by $P$ the Poisson bivector of the variety $O\times X_0$. Let
us show that  $P_x\in\bigwedge^2 T_xX_0\oplus \bigwedge^2 \g_*x$ for
 $x\in X_0$ in general position. We may assume that $X_0$ is affine and smooth. We shall
see now that for $f\in \K[X_0]\subset \K[X_0\times O], g\in
\K[O]\subset \K[X_0\times O]$ the equality $\{f,g\}=0$ holds. Note
that $f\in \K[X_0\times O]^G$ whence $f$ commutes with any
hamiltonian $H_\xi, \xi\in\g$. Note that $H_\xi$ is constant on
$X_0$. Since the space $\Span_\K(H_\xi,\xi\in\g)$ is $G$-stable,
$H_\xi\in\K[O]$. Moreover, $\K[O]$ is algebraic over the subalgebra
generated by $H_\xi,\xi\in\g$. It follows from the uniqueness
property of a lifting of a derivation that $\{\K[X_0],\K[O]\}=0$.

 Since $P$ is $G$-invariant, the projection of $P_x$  to $\bigwedge^2
 T_xX_0$  for  $x\in X_0$ depends only on the $X_0$-component of
 $x$. This projection is a Poisson bivector on $X_0$. On the other
 hand, this is the bivector corresponding to the bracket
 $\{\cdot,\cdot\}^{X_0}$. So $X_0$ is Poisson, and the Poisson structure
 on $O\times X_0$ is the product structure.

Now let $G$ be  not necessarily semisimple. It remains to show that
the action $Z(G)^\circ:X_0$ is  Hamiltonian with the moment map
$\mu:=\mu_{G,X}|_{X_0}-\eta$. Clearly, $\mu$ is
$Z(G)^\circ$-invariant. Recall that $G$ acts on $X_0\times O$ as on
the product of $G$-varieties and the action of $Z(G)^\circ$ on $O$
is trivial. Thus $v(H_\xi)_x=v(H_\xi|_{X_0})_x$ for $x\in X_0,
\xi\in \z(\g)$. In the LHS (resp., RHS) of the previous equality $v$
denotes the Hamiltonian vector field on $O\times X_0$ (resp.
$X_0\cong X_0\times \{eH^\circ\}$). Thus the functions
$H_\xi|_{X_0},\xi\in\z(\g),$ satisfy condition (H1). \end{proof}

There is the natural embedding $\z(\g)=\g^G\hookrightarrow \g\quo
G$.

\begin{corollary}\label{corollary:4.1.6}
Let $G,X$ be as above. Then
\begin{enumerate}\item $\im \psi_{G,X}\subset \z(\g)$ and
$\psi_{G,X}=\psi_{Z(G^\circ),X}$.\item
$\overline{\defe}_G(X)=\overline{\defe}_{Z(G)^\circ}(X)$,
$\underline{\defe}_G(X)=\underline{\defe}_{Z(G)^\circ}(X)$.
\end{enumerate}
\end{corollary}
\begin{proof}
The first assertion follows directly from the definitions. To prove
the second one we may assume that $X\cong X_0\times O$ (in the
notation of the previous theorem). In this case our assertions are
obvious. \end{proof}

\subsection{The affine case}\label{subsection_CN2}
We preserve the notation of the previous subsection and suppose that
$X$ is affine. Recall that $G$ is a connected group such that
$G\cong Z(G)^\circ\times (G,G)$.

We have a Hamiltonian open embedding $(X_0\times
O)/\Gamma\hookrightarrow X$ with the image $X^0$. Note that $X_0$ is
a closed subvariety in $X$. In particular, $X_0$ is  affine. Denote
by $\overline{O}$ the affine Hamiltonian variety $\Spec(\K[O])$.
$X_0\times O$ is embedded into $X_0\times \overline{O}$ as an open
subset with the complement of codimension not less than 2. Thus we
have the action $\Gamma:X_0\times \overline{O}$ by Hamiltonian
automorphisms and the Hamiltonian morphism $\iota:(X\times
\overline{O})/\Gamma\rightarrow X$ extending the embedding
$(X_0\times O)/\Gamma\hookrightarrow X$.

\begin{theorem}\label{theorem:4.2.1}
The morphism $\iota$ defined above is an isomorphism.
\end{theorem}
 In the proof of the theorem we use the following lemma proved, for
 example, in~\cite{Kraft}, Section 3.4.
 \begin{lemma}\label{lemma:4.2.2}
Let $X,Y$ be irreducible affine varieties and $\varphi:X\rightarrow
Y$ a birational morphism such that $\im\varphi$ contains an open
subset $Y^0\subset Y$ with $\codim_Y(Y\setminus Y^0)>1$. Suppose $Y$
is normal. Then $\varphi$ is an isomorphism.
\end{lemma}
\begin{proof}[of Theorem~\ref{theorem:4.2.1}]
It is enough to prove that $\codim_XX\setminus X^0\geqslant 2$
(Lemma~\ref{lemma:4.2.2}). In the proof we may replace $G$ with $(G,G)$
and assume that $G$ is semisimple. The proof is in three steps.

{\it Step 1.} Let us prove that $\mu_{G,X}|_{\overline{Gx}}$ is a
finite morphism for $x\in X^0$. This would imply, in particular,
that the closed $G$-orbit in $\overline{Gx}$ is a point for $x\in
X^0$ (and hence, in virtue of the Luna slice theorem, for any $x\in
X$).

Put $A= \K[\overline{Gx}], B= \K[\overline{G\eta}]$. The dominant
morphism $\mu_{G,X}:\overline{Gx}\rightarrow \overline{G\eta}$
induces the monomorphism $B\hookrightarrow A$. The corresponding
extension of the fraction fields $\Quot(B)\subset \Quot(A)$ is
finite, its degree is equal to $\#(G_\eta/G_x)$. Denote by
$\overline{A},\overline{B}$ the integral closures of $A$ and $B$ in
$\Quot(A)$. Since $\overline{A},\overline{B}$ are integrally closed
in $\Quot(A)$, one gets
$\Quot(\overline{A})=\Quot(\overline{B})=\Quot(A)$. Moreover, the
algebra extensions $A\subset \overline{A},B\subset \overline{B}$ are
finite and $\overline{A},\overline{B}$ are stable under the action
of $G$ on $\Quot(A)$. Put
$Z_1=\Spec(\overline{A}),Z_2=\Spec(\overline{B})$. These are normal
$G$-varieties. We obtain  the following commutative diagram

\begin{picture}(100,63)
\put(5,50){$Z_1$} \put(5,5){$\overline{Gx}$} \put(80,50){$Z_2$}
\put(80,5){$\overline{G\eta}$} \put(20,53){\vector(1,0){57}}
\put(40,55){\scriptsize$\varphi_1$}\put(20,8){\vector(1,0){57}}\put(40,10){\scriptsize$\varphi_2$}
\put(8,47){\vector(0,-1){30}}\put(84,47){\vector(0,-1){30}}
\put(10,28){\scriptsize$\psi_1$}\put(86,28){\scriptsize$\psi_2$}
\end{picture}

Here the morphisms  $\psi_1,\psi_2$ are finite and the morphisms
$\psi_1,\varphi_1$ are birational. Note that both $Z_1$ and $Z_2$
contain an open orbit isomorphic to $Gx$. Being birational and
$G$-equivariant, the morphism $\varphi_1$ induces an isomorphism of
these orbits. Further, note that the $G$-variety $\overline{G\eta}$
contains only finitely many orbits and the dimensions of all these
orbits are even.  Taking into account that $\psi_2$ is finite, we
get $\codim_{Z_2}(Z_2\setminus Gx)\geqslant 2$.
Lemma~\ref{lemma:4.2.2} implies that $\varphi_1$ is an isomorphism.
Thus $\varphi_2$ is finite.

{\it Step 2.} Here we prove that $\mu_{G,X}^{-1}(0)=X^{G}$. Clearly,
$\mu_{G,X}(X^G)=0$.  Let $x\in \mu_{G,X}^{-1}(0)$. At first, we show
that $\g_x^\perp\subset \g$ consists of nilpotent elements. Denote
by $X_1$ the stratum of $X$ (see
Proposition~\ref{proposition:2.4.1}) containing $x$. By
Example~\ref{example:3.2.4}, $X_1$ is a Hamiltonian subvariety. By
the definition of $X_1$, $x\in \overline{X_1}^{max}$.
Proposition~\ref{proposition:3.1.7.5} implies
\begin{equation}\label{eq:3.1:1.51}d_x\mu_{G,X_1}(T^P_xX_1)=\g_x^\perp.\end{equation}
Let $\mathcal{N}$ denote the cone in $\g$ consisting of all
nilpotent elements. Note that $\im\mu_{G,X_1}\subset \mathcal{N}$.
Thus  $\im d_x\mu_{G,X_1}$ coincides with the image of the morphism
of the tangent cones $T_xX_1\rightarrow \mathcal{N}$ induced by
$\mu_{G,X_1}$. So $\im d_x\mu_{G,X_1}\subset \mathcal{N}$. Using
(\ref{eq:3.1:1.51}), we see that $\g_x^\perp\subset \mathcal{N}$. By
Theorem 1 from~\cite{Bourbaki}, Chapter 7, $\S10$, $\g_x$ is a
parabolic subalgebra of $\g$. Since $G/G_x$ is  quasiaffine,  we get
$\g_x=\g$.

{\it Step 3.} Let us complete the proof. Denote by $Z$ an
irreducible component of $X\setminus X^0$. This is a $G$-stable
subvariety in $X$. Denote by $\eta_1$ an element from a unique open
$G$-orbit in $\overline{\mu_{G,X}(Z)}$. The intersection $Z\cap
X^{G}$ is non-empty because $Z\subset X$ is closed and any closed
$G$-orbit consists of one point (step  1). Thus $0\in \mu_{G,X}(Z)$.
Since the dimension of any  fiber of a morphism is not less that the
dimension of a general one,
\begin{equation}\label{eq:3.1:1.41}\dim Z\cap
\mu_{G,X}^{-1}(0)\geqslant \dim Z\cap\mu_{G,X}^{-1}(\eta_1)=\dim
Z-\dim G\eta_1.\end{equation} On the other hand,
\begin{equation}\label{eq:3.1:1.42}\dim Z\cap
\mu_{G,X}^{-1}(0)\leqslant \dim\mu_{G,X}^{-1}(0)=\dim
X^{G},\end{equation} thanks to step 2. By step 1,
\begin{equation}\label{eq:3.1:1.43}\dim X^{G}=\dim X\quo G\leqslant \dim X-m_G(X)=\dim
X-\dim G\eta.\end{equation} It follows from
(\ref{eq:3.1:1.41}),(\ref{eq:3.1:1.42}),(\ref{eq:3.1:1.43}) that
$\dim Z-\dim G\eta_1\leqslant \dim X-\dim G\eta$. Since the
dimension of any adjoint orbit is even,  $\dim G\eta\geqslant\dim
G\eta_1+2$. Therefore $\dim Z\leqslant \dim X-2$. \end{proof}

\begin{corollary}\label{corollary:4.2.3}
Let $G$ be an arbitrary connected reductive group and $X$ an
irreducible affine CN Hamiltonian variety. Then:
\begin{enumerate}
\item Any closed orbit for the action $(G,G):X$ is a
point. \item Any irreducible component of a fiber of $\pi_{(G,G),X}$
contains a dense $(G,G)$-orbit.
\item The restriction of $\pi_{(G,G),X}:X\rightarrow
X\quo (G,G)$ to $\mu^{-1}_{(G,G),X}(0)$ is a finite bijective
morphism. In particular, if $X$ is normal, then this is an
isomorphism.
\item $\underline{\defe}_G(X)=\underline{\defe}_{Z(G)^\circ}(X\quo
(G,G))$, $\overline{\defe}_G(X)=\overline{\defe}_{Z(G)^\circ}(X\quo
(G,G))$. Here $X\quo (G,G)$ is equipped with the structure of a
Hamiltonian $Z(G)^\circ$-variety according to
Example~\ref{example:3.2.10}.
\end{enumerate}
\end{corollary}
\begin{proof}
Let us check that $(G,G)x$ is closed iff $x\in X^{(G,G)}$ iff
$\mu_{(G,G),X}(x)=0$. Indeed, to prove this we may replace $X$ with its
normalization. Then we are done by the previous theorem.

Now to prove the third assertion  it is enough to show that the
morphism $\pi_{(G,G),X}|_{X^{(G,G)}}$ is finite. Let $\widetilde{X}$
be the normalization of $X$. Then the natural morphisms
$\widetilde{X}\quo (G,G)\rightarrow X\quo (G,G)$,
$\widetilde{X}^{(G,G)}\rightarrow X^{(G,G)}$ are finite and
dominant. To prove assertions 2-4  we may assume that $X$ is normal.
Also we may assume that $G$ is connected and $G\cong
Z(G)^\circ\times (G,G)$. Now our assertions are direct consequences
of Theorem~\ref{theorem:4.2.1} (in the notation of this theorem the
Hamiltonian $Z(G)^\circ$-variety $X\quo (G,G)$ is isomorphic to
$X_0/\Gamma$). \end{proof}

Now we consider the case of a smooth variety $X$.

\begin{proposition}\label{corollary:4.2.4}
Let $G=Z(G)^\circ\times (G,G)$, $X$ be smooth and affine and the action $G:X$
locally effective. In the preceding notation, $(G,G)\cong
\Sp(2m_1)\times\ldots\times\Sp(2m_k)$, $\overline{O}$ is the direct
sum of the tautological $\Sp(2m_k)$-modules, $X_0$ is smooth and the
action $\Gamma:X_0\times \overline{O}$ is free.
\end{proposition}
\begin{proof}
It follows from Theorem~\ref{theorem:4.1.1} that the morphism
$\pi_{\Gamma,X_0\times \overline{O}}:X_0\times
\overline{O}\rightarrow X$ is \'{e}tale in codimension 1. By
Proposition~\ref{proposition:3.3}, $\pi_{\Gamma,X_0\times
\overline{O}}$ is \'{e}tale. Therefore $X_0\times \overline{O}$ is
smooth and the action of $\Gamma$ is free. It follows from the Luna
slice theorem that $\overline{O}$ is a $(G,G)$-module. Note that $\overline{O}$
is symplectic as a $(G,G)$-module. It is enough to prove that if $G$
is simple, then $G\cong \Sp(2m)$ and $\overline{O}$ is the
tautological $\Sp(2m)$-module. In~\cite{Vinberg_connections} the
list of all linear representations of simple groups possessing a
dense orbit is given. Only one of these representations is
symplectic. \end{proof}

\subsection{An estimate on the dimension of a
fiber of $\psi_{G,X}$}\label{subsection_CN3}
\begin{proposition}\label{proposition:4.3.1}
Let $G$ be a connected reductive group, $X$  an irreducible CN
Hamiltonian $G$-variety.  Then the codimension of any fiber of
$\psi_{G,X}$ is not less than $\underline{\defe}_G(X)$.
\end{proposition}
\begin{proof}
 Using Corollary~\ref{corollary:4.1.6},
we may replace $G$ with $Z(G)^\circ$ and assume that $G=T$ is a torus.
Further, we may 
assume that $X$
is normal.

To prove the proposition in this case we need three lemmas

\begin{lemma}\label{lemma:4.3.2}
Let $T$ be a torus,  $X$ an irreducible affine $T$-variety and $Z$
an irreducible component of a fiber of $\pi_{T,X}$. If the action
$T:X$ is locally effective, then so is the action $T:Z$.
\end{lemma}
\begin{proof}[of Lemma~\ref{lemma:4.3.2}]
Assume the converse: there exists a non-trivial connected subgroup
$T_0\subset T$ acting trivially on $Z$. Choose a point $z\in Z$ not
lying in another component of $\pi_{T,X}^{-1}(y)$, where $
y=\pi_{T,X}(z)$. We see that
\begin{equation}\label{eq:4.3:1}\pi_{T_0,X}^{-1}(\pi_{T_0,X}(z))\subset
\pi_{T,X}^{-1}(y).\end{equation} But the action $T_0:Z$ is trivial.
This yields
\begin{equation}\label{eq:4.3:2}\pi_{T_0,X}^{-1}(\pi_{T_0,X}(z))\cap Z=\{z\}.\end{equation}
It follows from (\ref{eq:4.3:1}), (\ref{eq:4.3:2}) and the choice of
$z$ that $\pi_{T_0,X}^{-1}(\pi_{T_0,X}(z))=\{z\}$. This implies that
any fiber of the quotient morphism $\pi_{T_0,X}$ is trivial. In
other words, the action $T_0:X$ is trivial. Contradiction.
\end{proof}

\begin{lemma}\label{lemma:4.3.3}
Let $T$ be a torus and $X$ a smooth irreducible Hamiltonian
$T$-variety such that the action  $T:X$ is locally effective. Then
for all $\eta\in\im\mu_{T,X}$ any irreducible component of
$\mu_{T,X}^{-1}(\eta)$ contains a point $x$ such that
$d_x\mu_{T,X}:T_xX\rightarrow \t$ is surjective.
\end{lemma}
\begin{proof}[of Lemma~\ref{lemma:4.3.3}]
Let $x$ be a point lying in a  component $Y$ of
$\mu_{T,X}^{-1}(\eta)$ and not lying in another component. There
exists a $T$-stable open affine neighborhood of $x$ in $X$
(\cite{Sumihiro}). Thus one may assume that $X$ is affine.

By the choice of $x$, any irreducible component of
$\pi_{T,X}^{-1}(\pi_{T,X}(x))$ containing $x$ is contained in $Y$.
It follows from Lemma~\ref{lemma:4.3.2} that $Y$ contains a point
$x_1$ such that $\dim Tx_1=\dim T$.

Let $X'$ be a stratum (see Proposition~\ref{proposition:2.4.1}) of
$X$ containing $x_1$. By Example~\ref{example:3.2.4},
 $X'$ is a $T$-stable subvariety of $X$ and the action $T:X'$ is
 Hamiltonian with the
moment map $\mu_{T,X'}=\mu_{T,X}|_{X'}$.
Proposition~\ref{proposition:3.1.7.5} and the equality $\dim
Tx_1=\dim T$ imply that $\mu_{T,X'}=\mu_{T,X}|_{X'}$ is a submersion
in $x_1$. In particular, $\mu_{T,X}$ is a submersion in $x_1$.
\end{proof}

\begin{lemma}\label{lemma:4.3.4}
Let $T$ be a torus,  $X$  an irreducible affine  $T$-variety, $X_1$
a $T$-stable subvariety of $X$. Then $\dim X_1-m_T(X_1)\leqslant
\dim X-m_T(X)$.
\end{lemma}
\begin{proof}[of Lemma~\ref{lemma:4.3.4}]
We may assume that the action $T:X$ is locally effective.  Let $T_1$
denote the inefficiency kernel for the action $T:X_1$. Then  $\dim
X_1-m_T(X_1)=\dim X_1-\dim T+\dim T_1=\dim X-\dim T-(\dim X-\dim
T_1-\dim X_1)\leqslant \dim X-\dim T-(\dim X\quo T_1-\dim X_1)$. It
remains to notice that $\dim X_1=\dim X_1\quo T_1\leqslant \dim
X\quo T_1$. \end{proof}

The proof of the proposition  is by induction on $\dim X$. The case
$\dim X=0$ is obvious.

Let  $T_0$ denote the inefficiency kernel for the action $T:X$.
Using Remark~\ref{remark:3.1.3}, we may  replace $T$ with $T/T_0$ and
assume that the action $T:X$ is locally effective. Let us choose
$\alpha\in\t$ and prove that $\dim \mu_{T,X}^{-1}(\alpha)\leqslant
\dim X-\dim T$. Choosing a point on $\mu_{T,X}^{-1}(\alpha)$ and
replacing  $X$ with an invariant affine neighborhood of this point, we
may assume that $X$ is affine.

Let $X=X_0\cup X_1\cup\ldots\cup X_k$ be the stratification
introduced in Proposition~\ref{proposition:2.4.1}, where $X_0$ is an
open stratum. By the inductive assumption,
$\dim\mu_{G,X}^{-1}(\alpha)\cap X_i\leqslant \dim X_i-m_T(X_i)$ for
$i>0$. It follows from Lemma~\ref{lemma:4.3.4} that $\dim
X_i-m_T(X_i)\leqslant\dim X-\dim T$. It remains to show that $\dim
\mu_{T,X_0}^{-1}(\alpha)$ $\leqslant \dim X-\dim T$. By
Lemma~\ref{lemma:4.3.3}, for any irreducible component $Z$ of
$\mu_{T,X_0}^{-1}(\alpha)$ there exists $z\in Z$ such that
$\mu_{T,X_0}$ is a submersion in $z$. This implies $\dim Z=\dim
X-\dim T$. \end{proof}

\subsection{On the structure of $C_{G,X},\tau_{G,X}$}\label{subsection_CN4}
 Put $\a_{G,X}=\overline{\im\psi_{G,X}}$. The space
$\z(\g)\cong \g^G$ is naturally embedded into $\g\quo G$. Since $X$
is CN, $\a_{G,X}$ lies in $\z(\g)$ and coincides with
$\overline{\im\mu_{Z(G)^\circ,X}}=\a_{Z(G)^\circ,X}$. Denote by
$\tau^1_{G,X}$ the embedding $\a_{G,X}\rightarrow \g\quo G$ and by
$\tau^2_{G,X}:C_{G,X}\rightarrow \a_{G,X}$ a unique morphism such
that $\tau_{G,X}=\tau^1_{G,X}\circ\tau^2_{G,X}$.

\begin{proposition}\label{proposition:4.4.1}
Let $\t_0$ be the Lie algebra of the inefficiency kernel of the
action $Z(G):X$. Suppose $X$ is equidefectinal. Then $\a_{G,X}$ is
an affine subspace of $\z(\g)$ of dimension $\defe_G(X)$
intersecting $\t_0$ in a unique point.
\end{proposition}
\begin{proof}
The claim on the dimension of $\a_{G,X}$ is obvious. By
Corollary~\ref{corollary:4.1.6}, $\defe_G(X)=\defe_{Z(G)^\circ}(X)$.
Therefore we may replace $G$ with $Z(G)^\circ$ and assume that $G=T$ is
a torus. Denote by $T_0$ the connected subgroup of  $T$
corresponding to $\t_0$. Let $\xi$ be a point in
$\overline{\im\mu_{T_0,X}}$ and $Z$ be an irreducible component of
$\mu_{T_0,X}^{-1}(\xi)$. Since $T_0$ acts trivially on $X$,
$\mu_{T_0,X}^*(\K[\t_0])$ lies in the center of $\K(X)$. Thus $Z$ is
a component of a  Poisson subvariety of $X$.
Proposition~\ref{proposition:2.1.2} implies that $Z\subset X$ is a
Poisson subvariety. Since $T$ is connected, any Poisson subvariety
in $X$ is Hamiltonian (see Example~\ref{example:3.2.4}). For
$\xi\in\overline{\im\mu_{T_0,X}}$ in general position we get
$\underline{\defe}_T(Z)=m_T(Z)=\dim T-\dim T_0=\defe_T(X)$. Thus
$\defe_T(X)=\overline{\defe}_T(Z)+\dim
\overline{\im\mu_{T_0,X}}\geqslant
\defe_T(X)+\dim\overline{\im\mu_{T_0,X}}$. This implies that
$\im\mu_{T_0,X}$ is a point. Note that this point is the
(orthogonal) projection of $\im\mu_{T,X}$ to $\t_0$. Hence
$\im\mu_{T,X}$ is contained in an affine subspace in $\t$ of
dimension $\defe_T(X)$ intersecting $\t_0$ in a unique point.
Comparing the dimensions, we see that $\overline{\im\mu_{T,X}}$
coincides with this affine space. \end{proof}

The following proposition is the main result of this subsection.

\begin{proposition}\label{proposition:4.4.2}
Let $G$ be a connected reductive group, $X$  a normal irreducible
equidefectinal CN Hamiltonian $G$-variety. Then
$\im\widetilde{\psi}_{G,X}$ is an open subset of $C_{G,X}$ and the
restriction of $\tau^2_{G,X}:C_{G,X}\rightarrow \a_{G,X}$ to
$\im\widetilde{\psi}_{G,X}$ is \'{e}tale.
\end{proposition}
\begin{proof}
Recall that, by definition, $C_{G,X}$ is a normal variety. Since
$\overline{\defe}_G(X)=\underline{\defe}_G(X)$, the morphism
$\widetilde{\psi}_{G,X}$ is equidimensional by
Proposition~\ref{proposition:4.3.1}. An equidimensional morphism to
a normal variety is open (see~\cite{Chevalley}). In particular,
$\im\widetilde{\psi}_{G,X}$ is an open subset of $C_{G,X}$.

Since $G$ is connected, $\K(X)^G$ is algebraically closed in
$\K(X)$. Therefore $\K[C_{G,X}]$ is an integral closure of
$\psi^*_{G,X}(\K[\g]^G)$ in $\K(X)$. It follows from
Corollary~\ref{corollary:4.1.6} that
$\psi^*_{G,X}(\K[\g]^G)=\psi^*_{Z(G)^\circ,X}(\K[\z(\g)])$. In other
words, the varieties $C_{G,X}$ and $ C_{Z(G)^\circ,X}$ are naturally
isomorphic and
$\widetilde{\psi}_{G,X}=\widetilde{\psi}_{Z(G)^\circ,X},$
$\tau_{G,X}=\tau_{Z(G)^\circ,X}$. By the definitions of
$\a_{G,X},\tau^2_{G,X}$, we have $\a_{G,X}=\a_{Z(G)^\circ,X},$ $
\tau^2_{G,X}=\tau^2_{Z(G)^\circ,X}$. Therefore we may assume that
$G=T$ is a torus.

Denote by $T_0$ the inefficiency kernel for the action $T:X$. Since
$X$ is equidefectinal,  $\dim \a_{T,X}=\dim T/T_0$. Thus
$\psi^*_{T,X}(\K[\t])=\psi^*_{T/T_0,X}(\K[\t/\t_0])$, in other
words, $C_{T,X}$ and $C_{T/T_0,X}$ are naturally isomorphic. By
Proposition~\ref{proposition:4.4.1}, the restriction of the
projection $\t\rightarrow \t/\t_0$ to $\a_{T,X}$ is an isomorphism.
Clearly, we have the commutative diagram

\begin{picture}(120,65)
\put(5,53){$X$}\put(47,53){$C_{T,X}$}\put(100,53){$\a_{T,X}$}
\put(1,7){$C_{T/T_0,X}$}\put(75,7){$\t/\t_0\cong \a_{T/T_0,X}$}
\put(8,50){\vector(0,-1){32}}\put(107,50){\vector(0,-1){32}}
\put(48,50){\vector(-1,-1){33}}\put(15,57){\vector(1,0){30}}
\put(68,57){\vector(1,0){30}}\put(30,13){\vector(1,0){40}}
\put(109,35){ \scriptsize$\cong$}\put(41,35){\scriptsize $\cong$}
\end{picture}

Therefore we may assume that the action $T:X$ is (locally)
effective. Let us show now that $\tau_{T,X}:C_{T,X}\rightarrow \t$
is \'{e}tale in points of $C_{T,X}^{reg}\cap
\widetilde{\psi}_{T,X}(X^{reg})$.  Indeed, for any $y\in
C_{T,X}^{reg}$ and any $x\in X^{reg}\cap
\widetilde{\psi}_{T,X}^{-1}(y)$ we have
$d_x\mu_{T,X}=d_y\tau_{T,X}\circ d_x\widetilde{\psi}_{T,X}$. By
Lemma~\ref{lemma:4.3.3}, there exists a point $x\in X^{reg}\cap
\widetilde{\psi}_{T,X}^{-1}(y)$ such  that $d_x\mu_{T,X}$ is a
surjection. Thus $\tau_{G,X}$ is \'{e}tale in $y$.

Now we check that $C_{T,X}^{reg}\cap
\widetilde{\psi}_{T,X}(X^{reg})$ is an open subset in
$\im\widetilde{\psi}_{T,X}$, whose complement is of codimension not
less than 2. Indeed, $\codim_{C_{T,X}}C_{T,X}^{sing},$ $\codim_X
X^{sing} \geqslant 2$ because $X,C_{T,X}$ are normal. Since
$\widetilde{\psi}_{G,X}$ is equidimensional
(Proposition~\ref{proposition:4.3.1}), $\codim_{\im
\widetilde{\psi}_{T,X}}\im\widetilde{\psi}_{T,X}\setminus
\widetilde{\psi}_{T,X}(X^{reg})\geqslant 2$. This shows our claim.

To complete the proof of the proposition it is enough to apply
Proposition~\ref{proposition:3.3} to the morphism
$\im\widetilde{\psi}_{T,X}\rightarrow \t$. \end{proof}

\subsection{The proof of Theorem~\ref{theorem:1.2.9} for CN
varieties}\label{subsection_CN5} In this subsection we suppose that
$X$ is normal, affine and strongly equidefectinal.  In the proof of
theorem~\ref{theorem:1.2.9} we may assume that $G\cong
Z(G)^\circ\times (G,G)$.

\begin{proposition}\label{proposition:4.5.1}
Let $T$ be a torus and $X$ be  irreducible affine strongly
equidefectinal Hamiltonian $T$-variety. Then the action $T:X$ is
stable.
\end{proposition}
\begin{proof}
We may assume that the action $T:X$ is effective. Indeed, let $T_0$
be the inefficiency kernel. Since $X$ is equidefectinal, the
hamiltonians $H_\xi,\xi\in\t_0,$ are constant. Thus $X$ is strongly
equidefectinal also as a Hamiltonian $T/T_0$-variety (with the same
stratification as for the action $T:X$).

Let $X=\coprod_i X_i$ be a stratification given by
Definition~\ref{definition:1.2.8}. Let us show that there exists a
stratum $X_i$ such that $\defe_T(X_i)=\dim T$ and there is a closed
$T$-orbit in $X_i$. Indeed, otherwise all closed $T$-orbits lie in
$\bigcup_{i\in J}X_i$, where $J=\{i|
\defe_T(X_i)<\dim T\}$. In other words, $X\quo T=\bigcup_{i\in
J}\overline{X_i}\quo T$. Thus there exists $i\in J$ such that
$\overline{X_i}\quo T=\overline{X}\quo T$. But
$\dim\overline{\mu_{T,X}(\overline{X_i})}=\defe_T(X_i)<\dim\overline{\im\mu_{T,X}}$.
Therefore $\mu_{T,X}\quo T(\overline{X_i}\quo T)\neq \mu_{T,X}\quo
T(X\quo T)$. Contradiction. So there is a point $x\in X_i$, where
$\defe_T(X_i)=\dim T$, such that the orbit $Tx$ is closed.

The action $T:X$ is stable iff there is a closed orbit of dimension
$\dim T$, see~\cite{Popov}. It is enough to show that one can find
such an orbit even in $\overline{X_i}$. So we may assume that $x\in
X^{max}$.

Let us prove that the action of $T_0:=(T_x)^\circ$ on $T_xX/T^P_xX$
is trivial. Assume the converse. Let us choose a $T_0$-stable
complement $V$ to $T_x^PX$ in $T_xX$. It follows from the Luna slice
theorem that there is a $T_0$-stable smooth locally-closed
subvariety $Y\subset X$ such that $x\in Y, T_xY=V$. Replacing $Y$ with
some open subset we may assume that $T_y^PX\oplus T_yY=T_yX$ for any
$y\in Y$. By the choice of $Y$, $\xi_*y\in T_yY$ for any $\xi\in
\t_0$. But since the action $T_x:X$ is Hamiltonian, $\xi_*y\in
T^P_yX$ for any $\xi\in\t_0$. Thus the action $T_0:Y$ is trivial.

Applying the slice theorem again, we see that it is enough to show
that the action $T_0:T_xX/\t_*x$ is stable. Since $T_0$ acts
trivially on $T_xX/T^P_xX,\t_*x$ we reduce to the proof of the
stability of the action $T_0:T^P_xX$.

Let us prove that an action of a torus $T$ on a symplectic
$T$-module $U$ is stable. Indeed, we may assume that the action is
effective. Choose  linearly independent weights
$\lambda_1,\ldots,\lambda_k, k=\dim T,$ of the $T$-module $U$. Since
$U\cong U^*$, we see that $-\lambda_1,\ldots,-\lambda_k$ are also
weights of $U$. Choose nonzero weight vectors $v_{\lambda_1},\ldots,
v_{\lambda_k},v_{-\lambda_1},\ldots, v_{-\lambda_k}$ in the
corresponding weight subspaces. The orbit of  $\sum_i
(v_{\lambda_i}+v_{-\lambda_i})$ is closed and its dimension is equal
to $\dim T$.

Since $T^P_xX$ is a symplectic $T_0$-module, we are done.
\end{proof}

Let $X_0,\eta,\Gamma$ be such as in Subsection~\ref{subsection_CN1},
$\overline{O}$ such as in Subsection~\ref{subsection_CN2}.

\begin{lemma}\label{lemma:4.5.2}
The Hamiltonian action $Z(G)^\circ:X_0$ is strongly equidefectinal.
\end{lemma}
\begin{proof}
The Hamiltonian action $G:X_0\times O$ is strongly equidefectinal
because there is an \'{e}tale Hamiltonian morphism $X_0\times
O\rightarrow X$ (see Lemma~\ref{lemma:16.6},
Theorem~\ref{theorem:4.1.1}). Note that any $G$-stable locally
closed subvariety $Y\subset X_0\times O$ has the form $Y_0\times O$
for some locally closed $Z(G)^\circ$-stable subvariety $Y_0\subset
X_0$. Clearly, $Y$ is a Hamiltonian subvariety of $X$ iff $Y_0$ is a
Hamiltonian subvariety of $X_0$. If $Y$ is irreducible, then
$\underline{\defe}_G(Y)=\underline{\defe}_{Z(G)^\circ}(Y_0),
\overline{\defe}_G(Y)=\overline{\defe}_{Z(G)^\circ}(Y_0)$. The
intersections of the strata of $X$ with $X_0$ form a stratification
of $X_0$ satisfying the assumptions of
Definition~\ref{definition:1.2.8}. \end{proof}

\begin{proof}[of Theorem~\ref{theorem:1.2.9} in the CN case]
 By Theorem~\ref{theorem:4.2.1}, $X\cong (X_0\times
\overline{O})/\Gamma$. We easily reduce to the case $X=X_0\times
\overline{O}$. Clearly, $X$ satisfies condition (b) or (c) of the
second assertion iff $\eta=0$ iff $\overline{O}$ is a point. By
Lemma~\ref{lemma:4.5.2}, a Hamiltonian $Z(G)^\circ$-variety $X_0$ is
strongly equidefectinal. Proposition~\ref{proposition:4.5.1} implies
that the action $Z(G)^\circ:X_0$ is stable. This completes the proof
because $\overline{O}$ contains the dense $(G,G)$-orbit $O$.
 \end{proof}

\section{Reduction to the central-nilpotent case}
This is the most important section of the paper. Here we show how to
reduce the proofs of the theorems stated in the Introduction to the case
of CN varieties. Throughout the section $G$ denotes a reductive
group and $X$ a quasiprojective Hamiltonian $G$-variety, if
otherwise is not stated.

The first subsection is devoted to an algebraic version of the
local-cross section theorem of Guillemin and Sternberg,
see~\cite{GS}. The original theorem deals with Hamiltonian actions
of compact groups on real  manifolds. Knop in~\cite{Knop4} proved an
analog of this theorem for Hamiltonian actions of reductive
algebraic groups on symplectic varieties. Our approach is based on
his.

Let us explain what we mean by a cross-section. Suppose that $X$ is
quasi-projective and normal. Fix a Levi subgroup $L\subset G$. In
general, the subvariety $\mu_{G,X}^{-1}(\l)\subset X$ is
$N_G(L)$-stable but the action of $N_G(L)$ is not Hamiltonian.
However, there is an $N_G(L)$-stable open subset $\l^{pr}\subset\l$
 such that the $N_G(L)$-subscheme
$\widetilde{Y}=\mu_{G,X}^{-1}(\l^{pr})$ is normal
(Corollary~\ref{corollary:5.1.2}), has a natural Hamiltonian
structure with the moment map $\mu_{G,X}|_{\widetilde{Y}}$
(Proposition~\ref{App:4.13}) and   the natural morphism
$G*_{N_G(L)}\widetilde{Y}\rightarrow X$ is \'{e}tale
(Corollary~\ref{corollary:5.1.2}).

In Subsection~\ref{subsection_Red2} we use the construction of the
previous subsection for some special choice of $L$. Namely, take for
$L$ the centralizer of $\mu_{G,X}(x)_s$ for $x\in X$ in general
position. The variety $\widetilde{Y}=\mu_{G,X}^{-1}(\l^{pr})$ is
$N_G(L)$-irreducible and the natural morphism
$G*_{N_G(L)}\widetilde{Y}\rightarrow X$ is an open embedding
(Proposition~\ref{proposition:5.2.2}). The $N_G(L)$-Hamiltonian
variety $\widetilde{Y}$ is CN. Moreover, $\widetilde{Y}$ is affine
provided so is $X$. Choose an irreducible (=connected) component $Y$
of $\widetilde{Y}$. Put $\a_{G,X}^{(Y)}=\a_{N_G(L,Y),Y}$ (see
Subsection~\ref{subsection_CN4}) and $W_{G,X}^{(Y)}=N_G(L,Y)/L$. In
the case when $X$ is affine and equidefectinal these definitions
coincide with those given in the Introduction. At the end of
Subsection~\ref{subsection_Red2} we introduce the factorization of
the morphism $\tau_{G,X}=\tau^1_{G,X}\circ\tau^2_{G,X}$, see
the Introduction.

In Subsection 5.3. we prove Theorem~\ref{theorem:1.2.9}. The key
step is to reduce the proof for the pair $(G,X)$ to the proof for
the pair $(L,Y)$, where $L,Y$ are such as in the previous paragraph.

Unfortunately, this trick does not work for
Theorems~\ref{theorem:1.2.3},\ref{theorem:1.2.5}. For example, $Y$
does not intersect some fibers of $\psi_{G,X}$.  In
Subsections~\ref{subsection_Red4}-\ref{subsection_Red6} we construct
another CN Hamiltonian $L$-variety $R$ from $Y$.

In Subsection~\ref{subsection_Red4} we recall some (mostly standard)
properties of $H$-invariants of an affine $G$-varieties, where $H$
is the unipotent radical or the derived subgroup of a parabolic subgroup
of $G$.

To construct  $R$ we need some parabolic subgroup $P\subset G$ with
Levi subgroup $L$. If we want $R$ to have good properties, we should
make some special choice of $P$. In Subsection~\ref{subsection_Red5}
we establish the notion of a parabolic subgroup  $P$ compatible with
$Y$ and study the key property of such a subgroup
(Proposition~\ref{proposition:5.5.1}).

Subsection~\ref{subsection_Red6} is devoted to the construction of
$R$. A sketch of the construction was given in
Subsection~\ref{SUBSECTION_Intro3}.

In Subsection~\ref{subsection_Red7} the basic properties of $R$ are
studied. The next  subsection is devoted to the proofs of
Theorems~\ref{theorem:1.2.3},\ref{theorem:1.2.5}. Here we introduce
some good action of $W_{G,X}^{(Y)}$ on the quotient $R\quo L$
(Lemma~\ref{lemma:5.8.1}, Proposition~\ref{lemma:5.8.2}). Notice
that there is no natural action of $N_G(L,Y)$ on $R$.

In Subsection~\ref{subsection_Red10} we establish a relation between
$\K[C_{G,X}]$ and $\K[X]^G\cap \z(\K(X)^G)$
(Proposition~\ref{proposition:4.5.12}) and prove
Theorem~\ref{theorem:1.2.7}. In  Subsection~\ref{subsection_Red11}
we give an example of a connected group  $G$ and a conical
symplectic affine variety $X$ such that the group $W_{G,X}^{(Y)}$ is
not generated by reflections.

\subsection{Local cross-sections}\label{subsection_Red1}
Let $\l$ be a Levi subalgebra in $\g$. Put $N=N_G(\l)$.
\begin{definition}\label{definition:5.1.1}
An element $\xi\in\l$ is called {\it principal} if
$\z_\g(\xi_s)\subset\l$. The subset of  $\l$ consisting of all
principal elements is denoted by $\l^{pr}$.
\end{definition}

Clearly, $\l^{pr}$ is a nonempty open subset of $\l$. There is
another equivalent definition of $\l^{pr}\subset\l$. Namely,
consider a Cartan subalgebra $\t\subset\l$ and the root systems
$\Delta(\g),\Delta(\l)$ of $\g,\l$, respectively, associated with
$\t$. The subset $\l^{pr}\subset \l$ consists precisely of elements
$\xi\in\l$ such that the semisimple part $\xi_s$ of $\xi$ is
conjugate under the action of $N$ to some element of $\t^0$, where
\begin{equation}\label{eq:5.1:3}\t^0:=\t\setminus\bigcup_{\alpha\in
\Delta(\g)\setminus \Delta(\l)}\ker\alpha.\end{equation}

It follows from this definition and the Chevalley restriction
theorem that there is  $f\in \K[\l]^{N}$ such that $\l^{pr}=\{\xi\in
\l| f(\xi)\neq 0\}$.

Consider the subscheme
$\widetilde{Y}:=\mu_{G,X}^{-1}(\l^{pr})\subset X$. This is a
principal $N$-stable open subscheme in $\mu_{G,X}^{-1}(\l)$. In
particular, if $X$ is affine, then so is $\widetilde{Y}$.  Let us
describe some properties of the natural morphism
$\varphi:G*_N\widetilde{Y}\rightarrow X$.

\begin{proposition}\label{proposition:5.1.1}
 The morphism $\varphi:G*_N\widetilde{Y}\rightarrow X,
[g,y]\mapsto gy$, is  non-ramified. Its restriction to any
irreducible component of $G*_N\widetilde{Y}$ is dominant.
$\im\varphi=\psi_{G,X}^{-1}(Z)$ for some open subset $Z\subset\g\quo
G$.
\end{proposition}
\begin{proof}
The dimension of any irreducible component of $\mu_{G,X}^{-1}(\l)$
is not less than $\dim X-\dim\g+\dim\l$. Thus the dimension of any
component of $G*_N\widetilde{Y}$ is not less than $\dim X$.

Now we shall check that $\varphi$ is non-ramified, i.e. that for any
$y\in \widetilde{Y}$ the linear map
$d_y\varphi:T_y(G*_N\widetilde{Y})\rightarrow T_{\varphi(y)}X$ is
injective.

We have the natural identification $T_y(G*_N\widetilde{Y})\cong
T_y\widetilde{Y}\oplus \l^\perp_*y$. The restriction of $\varphi$ to
$\widetilde{Y}$ is the embedding $\widetilde{Y}\hookrightarrow X$.
Thus $d_y\varphi|_{T_y\widetilde{Y}}$ is an embedding. It remains to
show that $\xi\in\l$ provided $\xi_*y\in T_y\widetilde{Y},\xi\in
\g$. Let $T^P_yX$, $\omega_y$ be such as in
Subsection~\ref{subsection_Poisson2}. For any $\eta\in\l^\perp$ we
have
\begin{equation*}
0=\partial_{\xi_*y}H_{\eta}(y)=\omega_y(\xi_*y,\eta_*y)=\{H_\xi,H_\eta\}(y)=H_{[\xi,\eta]}(y)=([\mu_{G,X}(y),\xi],\eta).
\end{equation*}

Therefore $[\xi,\mu_{G,X}(y)]\in\l$. Let $\xi'$ denote the
projection of $\xi$ to $\l^\perp$. Recall that $\mu_{G,X}(y)\in
\l^{pr}$. Thus $\xi'\in\z_\g(\mu_{G,X}(y))\subset
\z_\g(\mu_{G,X}(y)_s)\subset \l$ whence $\xi\in\l$.

So $\varphi$ is non-ramified.  Comparing the dimensions, we see that
the restriction of $\varphi$ to any component of $G*_N\widetilde{Y}$
is dominant.

Now we prove the claim on the image of $\varphi$.  By the
alternative description of $\l^{pr}$ in the beginning of this
subsection,
\begin{equation}\label{eq:equivalence}x\in\im\varphi\Leftrightarrow G\mu_{G,X}(x)_s\cap\t^0\neq \varnothing.\end{equation}
 Put $Z=\t/W\setminus \pi_{W,\t}(\t\setminus \t^0)$, where $W=N_G(\t)/Z_G(\t)$.
Obviously,  $Z$ is an open subvariety in $\t/W\cong \g\quo G$. It
follows from (\ref{eq:equivalence}) that
$\im\varphi=\psi_{G,X}^{-1}(Z)$. \end{proof}

\begin{corollary}\label{corollary:5.1.2}
If $X$ is normal, then $\varphi:G*_N\widetilde{Y}\rightarrow X$ is
\'{e}tale and the scheme $\widetilde{Y}$ is normal.
\end{corollary}
\begin{proof}
$\varphi$ is \'{e}tale because this is a non-ramified dominant morphism
to a normal variety (see \cite{Miln}, Ch.1, Theorem 3.20). Thus the
scheme $G*_N\widetilde{Y}$ is a normal variety. By
Lemma~\ref{lemma:4.1.3}, $\widetilde{Y}$ is a normal scheme.
\end{proof}

We want to equip the $N$-variety $\widetilde{Y}$ with a natural
Hamiltonian structure.

 Let $y\in \widetilde{Y}$. We have the natural isomorphism
$T_yX\cong T_y\widetilde{Y}\oplus \l^\perp_*y$. The differential
$d_y\mu_{G,X}$ induces the isomorphism $\l^\perp_*y\cong
\l^\perp_*\mu_{G,X}(y)$. The restriction of the Kostant-Kirillov
form $\omega$ on the orbit $G\mu_{G,X}(y)$ (see
Example~\ref{example:2.3.1}) to the subspace
$V=\l^\perp_*\mu_{G,X}(y)\subset \g_*\mu_{G,X}(y)$ is nondegenerate.
Indeed, the subspaces $V,\l_*\mu_{G,X}(y)\subset \g_*\mu_{G,X}(y)$
are orthogonal with respect to this form. The form $\omega|_V$
induces the bivector $P^{\mu_{G,X}}_y\in \bigwedge^2 V$ (see
Subsection~\ref{subsection_Poisson2}).

Using the identification $\l^\perp_*y\cong V$, we obtain the
bivector $\widetilde{P}^{\mu_{G,X}}_y\in \bigwedge^2\l^\perp_*y$.
For $n\in N$ there is the natural isomorphism
$n_*:\l^\perp_*y\rightarrow \l^\perp_*(ny)$. Since
$\mu|_{\widetilde{Y}}$ is $N$-equivariant,
\begin{equation}\label{eq:bivector_inv}n_*\widetilde{P}^{\mu_{G,X}}_y=\widetilde{P}^{\mu_{G,X}}_{ny}.\end{equation}

For $f,g\in O_{X,y}$ define an element $\{f,g\}^Y\in
O_{\widetilde{Y},y}$ by
\begin{equation}\label{eq_red_bracket}
\{f,g\}^Y(y_1)=\{f,g\}(y_1)-\widetilde{P}^{\mu_{G,X}}_{y_1}(df\wedge
dg), y_1\in \widetilde{Y}.
\end{equation}

\begin{proposition}\label{App:4.13} Let $f,g\in O_{X,y}$. The element $\{f,g\}^Y\in O_{\widetilde{Y},y}$ depends only
on the restrictions of $f,g$ to $\widetilde{Y}$. Furthermore,
$\{\cdot,\cdot\}^Y:\K(\widetilde{Y})\otimes
\K(\widetilde{Y})\rightarrow \K(\widetilde{Y})$ is an $N$-invariant
Poisson bracket. The action $N:\widetilde{Y}$ is Hamiltonian with
the moment map $\mu_{N,\widetilde{Y}}=\mu_{G,X}|_{\widetilde{Y}}$.
\end{proposition}
\begin{proof}
Let $Y_0$ be an irreducible component of $\widetilde{Y}$. By
Proposition~\ref{proposition:5.1.1}, $\varphi(G*_N(NY_0))=GY_0$ is dense in $X$. Hence $Y_0\cap X^{max}\neq \varnothing$, so  we
may replace $X$ with $X^{max}$. Corollary~\ref{corollary:5.1.2}
implies that $G*_N\widetilde{Y}$ is smooth. Thus  $\widetilde{Y}$ is
smooth. Consider the distribution $T^PX$  and the  "2-form"
$\omega\in \operatorname{H}^0(X,\bigwedge^2 (T^{P}X)^*)$. By
Proposition~\ref{proposition:2.2.2},  $T^PX$ is involutory and
$\omega$ satisfies (\ref{eq_sympl_form}). Let $z\in \widetilde{Y}$.
 Choose  $\xi\in \l^\perp,
\eta\in T_z\widetilde{Y}\cap T^P_zX$. Since
$H_\xi|_{\widetilde{Y}}=0$, we have
$\omega_z(\xi_*z,\eta)=\partial_\eta H_\xi(z)=0.$ So $T^P_zX\cap
T_z\widetilde{Y}$ and $\l^\perp_*z$ are orthogonal with respect to
$\omega_z$. It follows that the bivector
$P^Y_z:=P_z-\widetilde{P}_z^{\mu_{G,X}}$ lies in $\bigwedge^2
T_z\widetilde{Y}$, equivalently, the bracket on $\widetilde{Y}$ is
well-defined. Note that
\begin{equation}\label{eq:5.1:5}T^{P^Y}_z\widetilde{Y}=T^P_zX\cap
T_z\widetilde{Y}, \g_*z\subset T^P_zX.\end{equation}

\begin{lemma}\label{lemma:5.1.5}
Let $H\subset G$ be algebraic groups, $Y$ a quasi-projective
$H$-variety and $\widetilde{V}$ a $G$-stable  distribution on
$G*_HY$ such that $\g_*[e,y]\in \widetilde{V}_{[e,y]}$ for all $y\in
Y$. Put $V_y=\widetilde{V}_{[e,y]}\cap T_yY$. If
 $\widetilde{V}$ is involutory, then so is $V$.
\end{lemma}
\begin{proof}[of Lemma~\ref{lemma:5.1.5}]
There is a quasi-section (see~\cite{VP}, $\S2$) $Z$ for the action
$H:G$ by the right translations such that the natural  morphism
$Z\rightarrow G/H$ is \'{e}tale. So $\varphi:Y\times Z\rightarrow G*_HY,
(y,z)\mapsto [z,y],$ is an \'{e}tale morphism. The distribution
$\varphi^*\widetilde{V}$ coincides with $V\otimes TZ\subset
T(Y\times Z)$.  \end{proof}

According to Lemma~\ref{lemma:5.1.5}, $T^{P^Y}_z\widetilde{Y}$ is an
involutory distribution. Let
$\omega^Y\in\operatorname{H}^0(\widetilde{Y},\bigwedge^2
(T^{P^Y}\widetilde{Y})^*)$ be the "2-form" corresponding to $P^Y$.
The form $\omega^Y$ is the restriction of $\omega$ to
$T^{P^Y}\widetilde{Y}$. Therefore $\omega^Y$ satisfies
(\ref{eq_sympl_form}). Applying Proposition~\ref{proposition:2.2.2},
we see that $P^Y$ is a Poisson bivector. Note that for $f\in
O_{X,y}$ the vector $v_y(f|_{\widetilde{Y}})$ coincides with the
projection of $v_y(f)$ to $T_y\widetilde{Y}$ with respect to the
decomposition $T_yX=T_y\widetilde{Y}\oplus \l^\perp_*y$. In
particular, for $\xi\in \l$ we have
$v(H_{\xi}|_{\widetilde{Y}})=\xi_*$ because $\widetilde{Y}$ is
$L$-stable.  Thus
$\mu_{G,X}|_{\widetilde{Y}}:\widetilde{Y}\rightarrow \l$ satisfies
the axioms of a moment map. \end{proof}

 The variety $\widetilde{Y}$ equipped
with this Hamiltonian action of $L$ is denoted by $\Red_G^N(X)$. Let
us investigate functorial properties of $\Red_G^N(\cdot)$.

\begin{proposition}\label{proposition:5.1.6}
Let $X_1,X_2$ be normal quasiprojective Hamiltonian $G$-varieties
and $\varphi:X_1\rightarrow X_2$ a Hamiltonian morphism. Then
$\varphi(\Red_G^N(X_1))\subset\Red_G^N(X_2)$. Let
$\Red_G^N(\varphi):\Red_G^N(X_1)\rightarrow \Red_G^N(X_2)$ be the
corresponding morphism. This is a Hamiltonian $N$-morphism.
\end{proposition}
\begin{proof}
The only not obvious thing here is that $\Red_G^N(\varphi)$ is a
Poisson morphism. Let $\{\cdot,\cdot\}^{Y_i}:\K(X_i)\otimes
\K(X_i)\rightarrow \K(\Red_G^N(X_i)),i=1,2,$ be the brackets defined
by (\ref{eq_red_bracket}). We have to prove that for $y\in
\Red_G^N(X_1)$ and $f,g\in O_{X_2,\varphi(y)}$ there is the equality
\begin{equation}\label{eq_functor}\{\varphi^*
f,\varphi^*g\}^{Y_1}(y)=\{f,g\}^{Y_2}(\varphi(y)).\end{equation}
Consider the bivectors $\widetilde{P}^{\mu_{G,X_1}}_y$ and
$\widetilde{P}^{\mu_{G,X_2}}_{\varphi(y)}$. Since $\varphi$ is a
Hamiltonian morphism, we see that
$d_y\Phi\widetilde{P}^{\mu_{G,X_1}}_y=P^{\mu_{G,X_2}}_{\varphi(y)}$.
 The equality (\ref{eq_functor}) follows from the definitions of $\{\cdot,\cdot\}^{Y_i}$,
 because $\varphi:X_1\rightarrow X_2$ is a Poisson morphism.
\end{proof}

Any irreducible component of $\Red_G^N(X)$ is a Hamiltonian
$L$-variety.

\begin{proposition}\label{proposition:5.1.7}
  $m_L(Y)=m_G(X)-\dim G+\dim
L$, $\overline{\defe}_G(X)=\overline{\defe}_L(Y)$,
$\underline{\defe}_G(X)=\underline{\defe}_L(Y)$ for any irreducible
component $Y\subset \widetilde{Y}$.
\end{proposition}
\begin{proof}
Since $\varphi:G*_{L}Y\rightarrow X$ is dominant and non-ramified,
$m_G(X)=\dim G/L+m_L(Y)$ and $GY$ is dense in $X$. The latter
implies the equality of the upper defects. Taking into account that
$\im\mu_{L,Y}\subset \l^{pr}$, we get $\rank_L(Y)=\rank_G(X)-\dim
G/N$. This implies
$\underline{\defe}_{G}(X)=\underline{\defe}_{L}(Y)$. \end{proof}

\subsection{$\a_{G,X}^{(Y)}$ and $W_{G,X}^{(Y)}$}\label{subsection_Red2}
We preserve the notation of the previous subsection. Suppose  that
$X$ is irreducible and normal. By Corollary~\ref{corollary:5.1.2},
if $\Red_G^N(X)$ is non-empty, then $\varphi$ is a dominant \'{e}tale
Hamiltonian morphism and $\Red_G^N(X)$ is a normal (possibly
non-connected) Hamiltonian $N$-variety. We want to show that for a
special choice of $L$ the morphism $\varphi$ is an open embedding.

Namely,  let $L\subset G$ be the centralizer in $G^\circ$ of
$\mu_{G,X}(x)_s$ for $x\in X$ in general position. In other words,
$L$ is the stabilizer of the closed orbit in general position for
the action $G^\circ:\overline{\im\mu_{G,X}}$ (so-called, principal
isotropy group, see~\cite{VP}, Theorem 7.12). Notice that $L$ is
defined uniquely up to $G^\circ$-conjugacy.

\begin{definition}\label{definition:5.2.1}
Such a subgroup $L$ is called the {\it principal centralizer} of
$X$.
\end{definition}

\begin{proposition}\label{proposition:5.2.2}
\begin{enumerate}
\item
The morphism $\varphi:G*_N\Red_G^N(X)\rightarrow X$ is an open
embedding. In particular,  $\Red_G^N(X)$ is $N$-irreducible. \item
If $X$ is, in addition, affine, then the natural morphism
$\Red_G^N(X)\quo N \rightarrow X\quo G$ is an open embedding.
\end{enumerate}
\end{proposition}
\begin{proof}
Let us check assertion 1.   Since $\varphi$ is \'{e}tale, it remains to
prove that $\varphi$ is injective. Let $y_1,y_2\in \Red_G^N(X),
g_1,g_2\in G$ be such that $g_1y_1=g_2y_2$. We may assume $g_1=1$.
Note that $\mu_{G,X}(y_1)_s=\Ad(g_2)\mu_{G,X}(y_2)_s$. The
centralizer of $\mu_{G,X}(y_i)_s$ in $\g$ coincides with $\l$. Thus
$g_2\in N$. Since $\varphi$ is an embedding, we see that
$G*_N\Red_G^N(X)$ is irreducible. By Lemma~\ref{lemma:4.1.3},
$\Red_G^N(X)$ is $N$-irreducible.

Recall that  $\im\varphi=\psi_{G,X}^{-1}(Z)$, where $Z$ is an open
subset of $\g\quo G$. Thus $G*_N \Red_G^N(X)\cong
\im\varphi=\pi_{G,X}^{-1}(\im\varphi\quo G)$. This implies assertion
3. \end{proof}

By the choice of $L$, for any $y\in \Red_G^N(X)$ the centralizer of
$\mu_{L,\Red_G^N(X)}(y)_s$ in $\g$ coincides with $\l$. Thus
$\Red_G^N(X)$ is a CN Hamiltonian $N$-variety.

Choose a component $Y\subset \Red_G^N(X)$ and put $N_0=N_G(L,Y)$. As
we have seen in the previous subsection, $Y$ is a normal CN
Hamiltonian $N_0$-variety with the moment map
$\mu_{N_0,Y}=\mu_{G,X}|_Y$. We say that $Y$ is an {\it
$L$-cross-section} of $X$. If $X$ is affine, then so is $Y$. It
follows from Proposition~\ref{proposition:5.2.2} that the triple
$(L,Y,N_0)$ is determined uniquely up to $G$-conjugacy and the
natural map $G*_{N_0}Y\rightarrow X$ is an open embedding, whose
image equals to $\psi_{G,X}^{-1}(Z)$ for some open subset
$Z\subset\g\quo G$.

Put $\a_{G,X}^{(Y)}=\a_{L,Y}$ (see Subsection~\ref{subsection_CN4}).
\begin{remark}\label{remark:5.2.3}
Suppose $X$ is equidefectinal. Then so is $Y$
(Proposition~\ref{proposition:5.1.7}). By
Proposition~\ref{proposition:4.4.1}, $\a_{G,X}^{(Y)}$ is an affine
subspace in $\z(\l)$. It follows directly from the definition of
$\a_{G,X}^{(Y)}$  that
$\a_{G,X}^{(Y)}=\a_{L,Y}=\overline{\pi(\mu_{G,X}(Y))}$, where $\pi$
is the projection $\l\rightarrow \l/[\l,\l]\cong \z(\l)$. So the
definition of $\a_{G,X}^{(Y)}$ given here coincides with that from
the Introduction.
\end{remark}

Put $W_{G,X}^{(Y)}=N_0/L$. Since $Y$ is a Hamiltonian $N_0$-variety,
$\a_{G,X}^{(Y)}\subset \z(\l)$ is stable under the natural action of
$W_{G,X}^{(Y)}$ on $\z(\l)$.

 Clearly, the restriction
of $f\in \K[\g]^G$ to $\a_{G,X}^{(Y)}\subset \g$ is
$W_{G,X}^{(Y)}$-invariant. So we have the natural morphism
$\tau^{1\,(Y)}_{G,X}:\a_{G,X}^{(Y)}/W_{G,X}^{(Y)}\rightarrow \g\quo
G$.

 The morphism
$\tau^{1\,(Y)}_{G,X}$ is, by definition, the composition of the
closed embedding $\a_{G,X}^{(Y)}/W_{G,X}^{(Y)}\hookrightarrow \l\quo
N_0$ and the  morphism $\l\quo N_0\rightarrow \g\quo G$ induced by
the restriction of functions. The last morphism is finite in virtue
of the Chevalley restriction theorem.  Thus $\tau^{1\,(Y)}_{G,X}$ is
finite.

Now we construct a $G$-invariant morphism
$\widehat{\psi}_{G,X}^{(Y)}:X\rightarrow
\a_{G,X}^{(Y)}/W_{G,X}^{(Y)}$ such that
$\psi_{G,X}=\tau^{1\,(Y)}_{G,X}\circ\widehat{\psi}^{(Y)}_{G,X}$. The
morphism $\psi_{N_0,Y}:Y\rightarrow \a_{G,X}^{(Y)}/W_{G,X}^{(Y)}$ is
$N_0$-equivariant. Thus we have a unique $G$-invariant morphism
$\psi:G*_{N_0}Y\rightarrow \a_{G,X}^{(Y)}/W_{G,X}^{(Y)}$ coinciding
with $\psi_{N_0,Y}$ on $Y$. By construction,
$\tau^{1\,(Y)}_{G,X}\circ\psi=\psi_{G,X}|_{G*_{N_0}Y}$. Since
$\tau^{1\,(Y)}_{G,X}$ is finite, the morphism $\psi$ can be extended
to the whole variety $X$. This extension is denoted by
$\widehat{\psi}_{G,X}^{(Y)}$.

\begin{remark}\label{remark:5.2.4}
The pair $(\a_{G,X}^{(Y)},W_{G,X}^{(Y)})$ depends on the choice of
$L,Y$ and so is determined uniquely up to $G$-conjugacy. However,
the quotient $\a_{G,X}^{(Y)}/W_{G,X}^{(Y)}$ and the morphisms
$\widehat{\psi}_{G,X}^{(Y)},\tau^{1\,(Y)}_{G,X}$ do not depend on
the choice of $L,Y$ in the following sense. Let $L'=gLg^{-1}, g\in
G^\circ,$ and $Y'$ be a component of $\Red_G^{N_G(L')}(X)$. There
exists an element $g_0$ such that
\begin{equation}\label{eq:5.2:1}
g_0Lg_0^{-1}=L', g_0Y=Y'.
\end{equation}
Moreover, $g_{01}^{-1}g_{02}\in N_0$ for any two elements $g_{01},g_{02}$
satisfying (\ref{eq:5.2:1}). Let $g_0\in N$ satisfy
(\ref{eq:5.2:1}). The isomorphism
$\iota:\a_{G,X}^{(Y)}/W_{G,X}^{(Y)}$ $\rightarrow
\a_{G,X}^{(Y')}/W_{G,X}^{(Y')}$  induced by $g_0$  does not depend
on the choice of $g_0$. Clearly,
$\iota\circ\widehat{\psi}^{(Y)}_{G,X}=\psi_{G,X}^{(Y')},
\tau_{G,X}^{1\,(Y)}=\tau_{G,X}^{1\,(Y')}\circ\iota$.

We write $\a_{G,X}/W_{G,X},\widehat{\psi}_{G,X},\tau^1_{G,X}$
instead of
$\a_{G,X}^{(Y)}/W_{G,X}^{(Y)},\widehat{\psi}_{G,X}^{(Y)},\tau_{G,X}^{1\,(Y)}$.
\end{remark}

By the definitions of $C_{G,X},\widetilde{\psi}_{G,X},\tau_{G,X}$
(see the Introduction) there is a unique morphism
$\tau^2_{G,X}:C_{G,X}\rightarrow \a_{G,X}/W_{G,X}$ such that
$\widehat{\psi}_{G,X}=\tau^2_{G,X}\circ\widetilde{\psi}_{G,X}$,
$\tau_{G,X}=\tau_{G,X}^1\circ\tau_{G,X}^2$.

In the case when $X$ is affine our morphisms are depicted on the
following commutative diagram
\begin{equation}\label{eq:5.2:2}
\begin{picture}(230,85)
\put(5,2){$X$}\put(5,40){$Y$}\put(1,70){$Y^{(L,L)}$}
\put(47,40){$Y\quo (L,L)$} \put(130,40){$Y\quo N_0$}
\put(50,2){$X\quo
G$}\put(190,40){$\a_{G,X}/W_{G,X}$}\put(195,70){$\a_{G,X}^{(Y)}$}\put(130,2){$C_{G,X}$}
\put(195,2){$\g\quo G$} \put(8,68){\vector(0,-1){20}}
\put(8,38){\vector(0,-1){27}} \put(12,6){\vector(1,0){38}}
\put(12,43){\vector(1,0){30}} \put(90,43){\vector(1,0){34}}
\put(60,37){\vector(0,-1){26}} \put(76,6){\vector(1,0){48}}
\put(160,43){\vector(1,0){28}}
\put(20,72){\vector(1,0){170}}\put(155,6){\vector(1,0){34}}
\put(200,67){\vector(0,-1){18}}\put(145,10){\vector(2,1){50}}
\put(130,38){\vector(-2,-1){55}}
\put(12,70){\vector(3,-2){30}}\put(200,37){\vector(0,-1){26}}
\put(30,59){\scriptsize $\cong$} \put(202,21){\scriptsize
$\tau^1_{G,X}$} \put(165,9){\scriptsize $\tau_{G,X}$}
\put(149,24){\scriptsize $\tau^2_{G,X}$} \put(60,78){\scriptsize
$\mu_{N_0/(L,L),Y^{(L,L)}}=\mu_{G,X}|_{Y^{(L,L)}}$}
\put(90,9){\scriptsize $\widetilde{\psi}_{G,X}\quo G$}
\end{picture}
\end{equation}

The isomorphism $Y^{(L,L)}\cong Y\quo (L,L)$ takes place by
Corollary~\ref{corollary:4.2.3}.

\subsection{The proof of Theorem~\ref{theorem:1.2.9}. The general
case.}\label{subsection_Red3}
\begin{lemma}\label{lemma:5.3.1}
Suppose $X$ is strongly equidefectinal. Let $L$ be a Levi subgroup
of $G$, $N=N_G(L)$. Then the Hamiltonian $N$-variety $\Red_G^N(X)$
is strongly equidefectinal.
\end{lemma}
\begin{proof}
First we suppose that $X=X^{max}$. We may assume that $X$ is
$G$-irreducible. $X$ is strongly equidefectinal iff $X$ is
equidefectinal. It follows from Corollary~\ref{corollary:5.1.2} that
$\Red_G^N(X)$ is smooth. By (\ref{eq:5.1:5}), the Poisson bivector
on any component of $\Red_G^N(X)$ has constant rank.
Proposition~\ref{proposition:5.1.7} implies that any component of
$\Red_G^N(X)$ is equidefectinal as a Hamiltonian $L$-variety. We are
done.

In the general case the assertion of the lemma follows from the fact
that $\Red_G^N(\cdot)$ is a functor, see
Proposition~\ref{proposition:5.1.6}. \end{proof}

We recall that a subset  $X^0\subset X$ is said to be $G$-saturated
if $X^0=\pi_{G,X}^{-1}(\pi_{G,X}(X^0))$.

\begin{proof}[of Theorem~\ref{theorem:1.2.9}] Let $L$ denote the principal centralizer of $X$, see
Definition~\ref{definition:5.2.1}. Choose an  $L$-cross-section
$Y\subset \Red_G^{N_G(L)}(X)$. The natural morphism
$G*_{N_G(L,Y)}Y\rightarrow X$ is an open embedding
(Proposition~\ref{proposition:5.2.2}) and its image is $G$-saturated
(Proposition~\ref{proposition:5.1.1}). Therefore we may assume that
$X=G*_{N_G(L,Y)}Y$. By Lemma~\ref{lemma:5.3.1}, $Y$ is a strongly
equidefectinal Hamiltonian $L$-variety. It remains to prove the
theorem for the action $L:Y$. Since $Y$ is a CN Hamiltonian
$L$-variety, we are done (see Subsection~\ref{subsection_CN5}).
\end{proof}

\subsection{$P_u$- and $(P,P)$-invariants}\label{subsection_Red4}
In this subsection  $X$ is an affine $G$-\!\! variety,  $G$ is a
connected reductive group. Let $P$ be a parabolic subgroup of $G$,
$L$ a Levi subgroup of $P$. Let $P_u$ denote the unipotent radical
of $P$. Put $P_0=(P,P)$. Recall that $P_0=P_u\leftthreetimes (L,L)$.

The algebra $ \K[G]^{P_u}$ is finitely generated (\cite{Grosshans}).
Thus $\K[G]^{P_0}$ is also finitely generated. Till the end of the
section $H$ denotes one of the groups $P_u,P_0$. Put
$\overline{G/H}=\Spec(\K[G]^{H})$.

The homogeneous space $G/H$ is quasiaffine since the character group
of $H$ is trivial. Fix an open $G$-equivariant embedding
$G/H\hookrightarrow \overline{G/H}$.

There is the isomorphism $(\K[G/H]\otimes
 \K[X])^G\cong  \K[X]^{H}$ induced by the restriction of
 functions from $(\K[G/H]\otimes \K[X])^G$ to $\{eH\}\times X\subset \overline{G/H}\times X$ (see~\cite{Popov}).
 Thus the algebra $\K[X]^{H}$ is finitely
 generated.

The subalgebras $\K[X]^{P_u}, \K[X]^{P_0}\subset \K[X]$ are stable
under the action of $L$ and  $(\K[X]^{P_u})^L=(\K[X]^{P_0})^L=
\K[X]^P= \K[X]^G$, $(\K[X]^{P_u})^{(L,L)}=\K[X]^{P_0}$.

Put $X\quo H=\Spec( \K[X]^{H})$. Let $\pi_{H,X}:X\rightarrow X\quo
H$ be the corresponding morphism.  Note that this morphism is
dominant but, in general, not surjective. We have a unique action
$L:X\quo H$ such that the morphism $\pi_{H,X}$ is
$P_u$-invariant and $L$-equivariant.

 The morphism
$\pi_{H,X}$ is the composition of the embedding $X=\{eH\}\times
X\hookrightarrow \overline{G/H}\times X$ and the quotient morphism
$\pi_{G,\overline{G/H}\times X}$. Thus if $X^0$ is an open
(respectively, closed) affine $G$-saturated subset of $X$, then
$X^0\quo H$ is identified with an open (respectively, closed) affine
subvariety in $X\quo H$ so that $\pi_{H,X^0}=\pi_{H,X}|_{X^0}$. Note
that the identification is $L$-equivariant.

\begin{remark}\label{lemma:5.4.1}
There are two natural actions of $L$ on $\overline{G/H}$. Firstly,
there is the restriction of the action $G:\overline{G/H}$ to $L$.
 Secondly, there is the
action $L:\overline{G/H}$ induced from the action  $L:G/H$ by the
right translations.  This action commutes with the action of $G$ and
hence induces the action $L:X\quo H\cong (X\times
\overline{G/H})\quo G$ considered above.

If otherwise is not stated, we consider the  action of the first
type. Note, however, that the $L$-orbit of $eH$ is the same for the
both actions.
\end{remark}

In Subsection~\ref{subsection_Red5} we need to know whether there
exists
 $\lim_{t\rightarrow 0}\tau(t)eP_0$ in
$\overline{G/P_0}$, where $\tau$ is a one-parameter subgroup,
$\tau:\K^\times\rightarrow Z(L)$.

Fix a Cartan subalgebra $\t\subset\l$ and the corresponding root
system $\Delta(\g)$.
\begin{definition}\label{definition_compatible_roots} A system of simple
roots $\alpha_1,\ldots,\alpha_r\in\Delta(\g)$ is said to be {\it
compatible} with  $P$, if the inclusion $\g^\alpha\subset \p$ is
equivalent to $\alpha=\sum_{i=1}^r n_i\alpha_i,$ $
n_{l+1},\ldots,n_r\geqslant 0$.\end{definition} Fix a system of
simple roots  $\alpha_1,\ldots\alpha_r$ compatible with $P$ and let
$\pi_1,\ldots,\pi_r$ be the corresponding system of fundamental
weights.

\begin{lemma}\label{lemma:5.4.2}
Let $\tau:\K^\times\rightarrow Z(L)$ be a one-parameter subgroup.
Put $\xi:=\frac{d}{dt}\tau(t)|_{t=0}\in\z(\l)$. Then the limit
$\lim_{t\rightarrow 0}\tau(t)P_0$ exists iff $\xi\in [\g,\g]$ and
for all $i>l$ the inequality $\langle \pi_i,\xi \rangle\geqslant 0$
holds.
\end{lemma}
\begin{proof}
$\im\tau\subset(G,G)$ because the limit $\lim_{t\rightarrow
0}\pi(\tau(t))$ exists in $G/(G,G)$,
 where $\pi$ denotes  the projection $G\rightarrow G/(G,G)$.
Thus we may assume that $G$ is semisimple. Replacing $G$ with a
covering and $\tau$ with a positive multiple, we may assume that $G$
is simply connected. Let $V_i$ be the irreducible $G$-module with
the highest weight  $\pi_i$, and $v_i\in V_i$ be a highest vector,
$i=\overline{1,r}$. Put $v=v_{l+1}+\ldots +v_r\in
V_{l+1}\oplus\ldots\oplus V_r$. There is a unique $G$-equivariant
morphism $\overline{G/P_0}\rightarrow \overline{Gv}$ such that
$eP_0\mapsto v$. This is an isomorphism, see~\cite{VP_var}, Theorem
6. The limit $\lim_{t\rightarrow 0}\tau(t)v$ exists iff  $\langle
\xi,\pi_i\rangle\geqslant 0$ for all $i>l$. \end{proof}

\subsection{Compatible parabolic subgroups}\label{subsection_Red5}
Let $X$ be a normal irreducible affine Hamiltonian $G$-variety, $L$
the principal centralizer and $Y$ an $L$-cross-section  of $X$,
$N_0=N_G(L,Y)$, $L_0$ the connected component of the inefficiency
kernel of the action $L:Y\quo (L,L)$.

We need to choose some parabolic subgroup $P\subset G$ {\it
compatible} with $Y$. Fix a Cartan subalgebra $\t\subset \l$ and the
corresponding root system $\Delta(\g)$.

\begin{specdefi} Let us
embed the $L/(L,L)$-variety $Y^ {(L,L)}$ into some $L/(L,L)$-module
$V$. Choose a point $y_0\in Y^{(L,L)}$ with $(L_{y_0})^\circ=L_0$
(we recall that the $L$-varieties $Y^{(L,L)}$ and $Y\quo (L,L)$ are
isomorphic, see Corollary~\ref{corollary:4.2.3}). The dimension of
the support $S_{y_0}$ of $y_0$ (i.e. the convex hull of
$L/(L,L)$-weights of $y_0\in V$) equals $\rank L-\rank
L_0=\defe_G(X)$ (the last equality follows from
Corollary~\ref{corollary:4.2.3}). To any point $\zeta\in S_{y_0}$ we
assign a unique face $C_\zeta$ of a Weyl chamber of the dual root
system $\Delta^\vee(\g)$ such that $\zeta$ is contained in the
interior of $C_\zeta$.  Fix a point $\zeta\in S_{y_0}$ such that
$C_\zeta$ is maximal with respect to the inclusion among all
$C_{\zeta'}, \zeta'\in S_{y_0}$. Put
$$\q=\t\oplus \bigoplus_{\alpha,\langle
\alpha^\vee,\zeta\rangle\leqslant 0} \g^\alpha.$$  Note that $\l$ is
identified with a Levi subalgebra in $\q/\q_u$ and that for $\alpha\in\Delta(\g)$
\begin{equation}\label{eq:5.5:1}
 \langle \alpha^\vee,\zeta
\rangle=0\Rightarrow\alpha^\vee\in \l_0.
\end{equation}
Indeed, if $\alpha^\vee\not\in\l_0$, then $\alpha^\vee y_0\neq 0$.
In other words,  $S_{y_0}\not\subset \ker\alpha^\vee$. By the choice
of $\zeta$, $\langle\alpha^\vee,\zeta\rangle\neq 0$.

Choose a parabolic subalgebra $\underline{\p}\subset \q/\q_u$ such that $\l$
is a Levi subalgebra of $\underline{\p}$.  Let $\p$ be the inverse image of
$\p_0$ in $\q$ under the projection $\q\rightarrow \q/\q_u$. Let
$P,Q$ be the parabolic subgroups of $G$ corresponding to the
subalgebras $\p,\q\subset \g$. Such a parabolic subgroup $P\subset
G$ is said to be {\it compatible} with $Y$. Note that $P$ depends on
the choices of an $L/(L,L)$-module $V$, an embedding
$Y^{(L,L)}\hookrightarrow V$, a point $y_0$, an element $\zeta$ and
a subalgebra $\underline{\p}\subset \q/\q_u$.
\end{specdefi}

\begin{proposition}\label{proposition:5.5.1}
Let $P$ be a parabolic subgroup of $G$ compatible with $Y$. The
restriction of $\pi_{P_0,X}$ to $Y^{(L,L)}$ is generically finite.
Moreover, for some open subset $Y^1\subset Y^{(L,L)}$ the following
condition is satisfied:

if $y_1\neq y_2\in Y^1$ and $\pi_{P_0,X}(y_1)=\pi_{P_0,X}(y_2)$,
then there exists $g\in (N_0\cap G^\circ)\setminus L$ such that
$y_1=gy_2$.
\end{proposition}

The proof is based on the following lemma.

\begin{lemma}\label{lemma:5.5.2}
Let $y_0,\zeta,Q$ be such as in the previous
construction-definition, $\widetilde{L}=Z_{G^\circ}(\zeta), Q_0=(Q,Q)$. Then
\begin{enumerate}
\item $(\widetilde{L},\widetilde{L})\cap L\subset L_0$.
\item $\dim L(y_0,eQ_0)=\dim L-\dim L\cap
(\widetilde{L},\widetilde{L})$.
\item $L(y_0,eQ_0)$ is closed in $Y^{(L,L)}\times
\overline{G/Q_0}$.
\end{enumerate}
\end{lemma}
\begin{proof}
Let us check that $L\cap (\widetilde{L},\widetilde{L})\subset L_0$.
The derived subgroups of these two groups coincide. The space $\t\cap
[\widetilde{\l},\widetilde{\l}]$ is spanned by $\alpha^\vee$,
$\alpha\in\Delta(\widetilde{\l})$.  By (\ref{eq:5.5:1}), $\t\cap
[\widetilde{\l},\widetilde{\l}]\subset \t\cap\l_0$ whence
$\l\cap[\widetilde{\l},\widetilde{\l}]\subset\l_0$.

Note that $y_0$ is $L_0$-invariant and that $L_{eQ_0}=L\cap
(\widetilde{L},\widetilde{L})$. The equality $\dim L(y_0,eQ_0)=\dim
L-\dim L\cap (\widetilde{L},\widetilde{L})$ follows from assertion
1.

Let $\tau:\K^\times\rightarrow Z(\widetilde{L})$ be a one-parameter
subgroup such that the limit $\lim_{t\rightarrow
0}\tau(t)(y_0,eQ_0)$ exists. Choose a system of simple roots
$\alpha_1,\ldots,\alpha_r\in \Delta(\g)$ compatible with  $Q$ (see
Definition~\ref{definition_compatible_roots}) and the corresponding
system  $\pi_1,\ldots,\pi_r$  of fundamental weights. Put $l=\rank
[\widetilde{\l},\widetilde{\l}]$. It follows from the construction
of $Q$ that  $\zeta=\sum_{i>l} a_i\pi_i$ with $a_i<0$. In other
words, for any $\xi\in [\g,\g]\cap \t$ the inequalities
\begin{equation}\label{eq:7.11}\langle \pi_i,\xi\rangle\geqslant 0,\forall
i>l,
\end{equation} and $\langle \xi,\zeta \rangle\geqslant 0$ imply $\xi\in
[\widetilde{\l},\widetilde{\l}]$.

 By Lemma~\ref{lemma:5.4.2},
 $\xi:=\frac{d}{dt}|_{t=0}\tau$ lies in $ [\g,\g]$ and satisfies (\ref{eq:7.11}).
 Since the limit
$\lim_{t\rightarrow 0}ty_0$ exists,  $\langle
\xi,\zeta\rangle\geqslant 0$ (let us recall that $\zeta$ is
contained in the support of $y_0$). Thus $\xi\in
[\widetilde{\l},\widetilde{\l}]$. Since
$\xi=\frac{d}{dt}|_{t=0}\tau\in\z(\widetilde{\l})$,  we get $\xi=0$.

The Hilbert-Mumford theorem implies that the
$Z(\widetilde{L})$-orbit of $(y_0,eQ_0)$ is closed in
$Y^{(L,L)}\times \overline{G/Q_0}$. Since
$(\widetilde{L},\widetilde{L})\cap L$ leaves $(y_0,eQ_0)$ invariant,
the $L$-orbit of $(y_0,eQ_0)$ coincides with the
$Z(\widetilde{L})$-orbit. This proves assertion 3. \end{proof}

\begin{proof}[of Proposition~\ref{proposition:5.5.1}]
Without loss of generality we may assume that $G$ is connected.  The
subset $G*_{N_0}Y\subset X$ is affine, open and $G$-saturated.
Therefore we reduce to the case $X=G*_{N_0}Y$. In this case
$X\times\overline{G/P_0}=(G*_{N_0}Y)\times \overline{G/P_0}\cong
G*_{N_0}(Y\times \overline{G/P_0})$ (the last isomorphism is given
by $([g,y],z)\mapsto [g,(y,g^{-1}z)], g\in G,y\in
Y,z\in\overline{G/P_0}$). Clearly, $$\pi_{N_0,Y\times
\overline{G/P_0}}=\pi_{N_0/L, (Y\times \overline{G/P_0})\quo L}\circ
\pi_{L,Y\times \overline{G/P_0}}.$$  Hence it is enough to prove
that the restriction of $\pi_{L,Y\times \overline{G/P_0}}$ to
$Y^{(L,L)}\times \{eP_0\}$ is generically injective. Since
$Y^{(L,L)}\subset Y$ is a closed $L$-stable subvariety, it is enough
to show the analogous claim for the morphism $\pi_{L,Y^{(L,L)}\times
\overline{G/P_0}}$.

Let  $Q,y_0$ be such as in Construction-definition above. There is a
unique $G$-equivariant morphism $\overline{G/P_0}\rightarrow
\overline{G/Q_0}$ such that $gP_0\mapsto gQ_0$ for all $g\in G$. It
is enough to check that the restriction of $\pi_{L,Y^{(L,L)}\times
\overline{G/Q_0}}$ to $Y^{(L,L)}\times \{eQ_0\}$ is generically
injective. It is so provided the following two claims take place:

1) $ L_{(y_0,eQ_0)}$ is the inefficiency kernel for the action
$L:\overline{L(Y^{(L,L)}\times eQ_0)}$.

2) The orbit $L(y_0,eQ_0)$ is closed.

The second claim is assertion 3 of Lemma~\ref{lemma:5.5.2}. By
assertion 1 of Lemma \ref{lemma:5.5.2}, any point of
$Y^{(L,L)}\times eQ_0$ (and thence of $\overline{L(Y^{(L,L)}\times
eQ_0)}$)  is $L\cap (\widetilde{L},\widetilde{L})$-invariant.
Assertion 2 of Lemma~\ref{lemma:5.5.2} implies the first claim.
\end{proof}

\subsection{$P_u$-reduction. The construction}\label{subsection_Red6}
We use the notation introduced in the beginning of the previous
subsection.

\begin{lemma}\label{lemma:5.6.1}
Let $\p$ be a parabolic subalgebra in $\g$ with a Levi subalgebra
$\l$. Then the subvariety $\overline{P_uY}\subset X$ is an
irreducible component of $\mu_{G,X}^{-1}(\p)$.
\end{lemma}
\begin{proof}
 $\mu_{G,X}^{-1}(\l^{pr}+\p_u)$ is an open subvariety
of $\mu_{G,X}^{-1}(\p)$. It follows from the definition of $\l^{pr}$
that $\z_\g(\xi)\cap \p_u=\{0\}$ for any $\xi\in\l^{pr}$. Therefore
  $\l^{pr}+\p_u=P_u\l^{pr}$ and the action $P_u:\l^{pr}+\p_u$ is free. We deduce that the morphism $P_u\times \l^{pr}
  \rightarrow\l^{pr}+\p_u, (g,y)\mapsto gy,$ is an isomorphism. This implies that the
  morphism of schemes
  $P_u\times\mu_{G,X}^{-1}(\l^{pr})\rightarrow
\mu_{G,X}^{-1}(\l^{pr}+\p_u), (g,y)\mapsto gy,$ is an isomorphism
too. Since $Y$ is a connected (=irreducible) component of
$\mu_{G,X}^{-1}(\l^{pr})$, the subset $P_uY\subset X$ is open in
some component of $\mu_{G,X}^{-1}(\p)$. \end{proof}

Fix a parabolic subgroup $P\subset G$ compatible with $Y$.

Let $A$ denote the subalgebra of $ \K[X]$ generated by $
\K[X]^{P_u}$ and $\{H_\xi, \xi\in\l\}$. This is a finitely generated
$L$-stable subalgebra in $\K[X]$. Denote by
$\widetilde{\pi}_{P_u,X}$  the morphism $X\rightarrow \Spec(A)$
induced by the embedding $A\hookrightarrow \K[X]$, by $Z$ an
irreducible component of $\mu_{G,X}^{-1}(\p)$ and by $I_Z$ the ideal
of functions from $ \K[X]$ vanishing on $Z$. Note that $I_Z$ is an
$L$-stable ideal in $\K[X]$. Put $A_Z=A/(A\cap I_Z)$. This is a
subalgebra of $\K[Z]$ consisting of the restrictions of elements
from $A$ to $Z$. There is the natural action of $L$ on $A_Z$.

\begin{lemma}\label{lemma:5.6.2}
\begin{enumerate}
\item
For any $f,g\in A$ the restriction of $\{f,g\}|_Z$ depends only on
$f|_Z,g|_Z$ and is contained in $A_Z$. So $A_Z$ becomes a Poisson
algebra.
\item $L:\Spec(A_Z)$ is a Hamiltonian action with the hamiltonians $H_\xi|_Z, \xi\in\l$.
\end{enumerate}
\end{lemma}
\begin{proof}
Firstly, we check that $\{A,I_Z\}\subset I_Z$. Note that $I_Z$ is a
minimal prime ideal of the ideal $I=\Span_{\K[X]}(H_\xi,
\xi\in\p_u)$. Applying Lemma~\ref{lemma:2.1.4} to the algebra
$\K[X]/I$ and the ideal $I_Z/I$, we see that it is enough to show
that $\{A,I\}\subset I$.  If $\eta\in \p_u$, then
$\{H_\xi,H_\eta\}=H_{[\xi,\eta]}\in I$ and $\{f,H_\eta\}=-\eta_*f=0$
for $f\in \K[X]^{P_u}$. Since $\K[X]^{P_u},H_\xi,\xi\in\l$, generate
$A$, we get $\{A,I\}\subset I$.

To prove the first assertion of the lemma it is enough to show that
$A$ is a Poisson subalgebra of $\K[X]$. To do this we have to check
that the brackets of generators of $A$ lie in $A$. Let $f,g\in
\K[X]^{P_u}, \xi,\eta\in \l$. One checks directly that
$\{f,g\},\{H_\xi,f\}\in \K[X]^{P_u},
\{H_\xi,H_\eta\}=H_{[\xi,\eta]}$.

Assertion 2 is verified directly using
Definition~\ref{definition:3.1.1}. \end{proof}

 By
Lemma~\ref{lemma:5.6.1}, we may apply the previous construction to
$Z=\overline{P_uY}$.

\begin{definition}\label{definition:5.6.3}
By the {\it $P_u$-reduction} of $X$ associated with $Y$ we mean the
normalization of the Hamiltonian $L$-variety
$\Spec(A_{\overline{P_uY}})$.
\end{definition}

Till the end of the section $R$ denotes  the $P_u$-reduction of $X$
associated with $Y$ and $Z$ denotes the subvariety
$\overline{P_uY}\subset X$. $R$ is equipped with the natural
structure of a Hamiltonian $L$-variety, see
Example~\ref{example:3.2.3}. To make the notation less bulky, we
write $\underline{R}_Z$ instead $\Spec(A_Z)$ and $L'$ instead of
$(L,L)$.

\subsection{$P_u$-reduction. The basic properties}\label{subsection_Red7}
We preserve the notation of the previous subsection.

Let us make some remarks on morphisms between our varieties.
Firstly, we have the natural dominant morphism
$\widetilde{\pi}_{P_u,X}|_{Z}:Z\rightarrow \underline{R}_Z$. This
morphism is $P_u$-invariant and $L$-equivariant. The restriction of
this morphism to $Y$ is dominant and $L$-equivariant. Since $Y$ is
normal, this restriction can be lifted to an $L$-equivariant
dominant morphism $\widehat{\pi}:Y\rightarrow R$.

Secondly, we have the $L$-equivariant morphism $\nu:R\rightarrow
X\quo P_u$ corresponding to the composition of the homomorphisms
$\K[X]^{P_u}\twoheadrightarrow \K[X]^{P_u}/(I_Z\cap
\K[X]^{P_u})\hookrightarrow A_Z\hookrightarrow \K[R]$.

\begin{lemma}\label{lemma:5.7.1} The following diagram is
commutative. Here all horizontal arrows are quotient morphisms, the
morphism $Y\rightarrow R$ is $\widehat{\pi}$,  $R\rightarrow
\underline{R}_Z$ is the normalization, the morphism $R\rightarrow
X\quo P_u$ coincides with $\nu$, the morphisms $Y\rightarrow
Z\rightarrow X$ are embeddings, all vertical arrows in the rectangle
with the vertices $Y\quo L', X\quo P_0,X\quo P,Y\quo L$ are
determined uniquely by the commutativity condition, the morphism
$Z\rightarrow \underline{R}_Z$ is induced by the embedding
$A_Z\hookrightarrow \K[Z]$, all morphisms to $\l\quo L$ are of the
form $\psi_{L,\bullet}\quo L$, the morphism $X\quo G\rightarrow
\g\quo G$ is $\mu_{G,X}\quo G$ and $\l\quo L\rightarrow \g\quo G$
is induced by the restriction of functions.

\begin{equation}\label{eq:5.7:1}
\begin{picture}(250,140)
\put(30,130){$Y$}\put(80,130){$Y\quo L'$}\put(160,130){$Y\quo
L$}\put(230,110){$\l\quo L$} \put(30,90){$R$}\put(80,90){$R\quo
L'$}\put(160,90){$R\quo
L$}\put(30,50){$\underline{R}_Z$}\put(78,50){$\underline{R}_Z\quo
L'$}\put(158,50){$\underline{R}_Z\quo L$}\put(2,70){$Z$}
\put(2,5){$X$}\put(28,5){$X\quo P_u$}\put(82,5){$X\quo P_0$}
\put(140,5){$X\quo P\cong X\quo G^\circ$}\put(230,5){$X\quo G$}
\put(230,70){$\g\quo G$} \put(30,128){\vector(-1,-2){23}}
\put(5,68){\vector(0,-1){53}} \put(33,88){\vector(0,-1){28}}
\put(33,128){\vector(0,-1){28}}\put(33,46){\vector(0,-1){30}}
\put(93,88){\vector(0,-1){28}}
\put(93,128){\vector(0,-1){28}}\put(93,46){\vector(0,-1){30}}
\put(170,88){\vector(0,-1){28}}
\put(170,128){\vector(0,-1){28}}\put(170,46){\vector(0,-1){30}}
\put(40,134){\vector(1,0){36}}\put(40,94){\vector(1,0){36}}
\put(40,55){\vector(1,0){36}}
\put(110,134){\vector(1,0){44}}\put(110,94){\vector(1,0){44}}
\put(110,55){\vector(1,0){44}} \put(10,8){\vector(1,0){17}}
\put(55,8){\vector(1,0){25}} \put(110,8){\vector(1,0){27}}
\put(205,8){\vector(1,0){23}}\put(240,16){\vector(0,1){50}}
\put(240,108){\vector(0,-1){25}}\put(184,134){\vector(3,-1){42}}
\put(184,95){\vector(3,1){42}} \put(186.5,58){\vector(1,1){46}}
\put(10,73){\vector(1,-1){16}}
\end{picture}
\end{equation}
\end{lemma}

\begin{proof}
The commutativity of the piece of the diagram  inside the triangle
with vertices $Y\quo L, \l\quo L,$ $\underline{R}_Z\quo L$ stems
from the constructions of the moment maps for the actions of $L$ on
$Y,\underline{R}_Z,R$ (see Subsections~\ref{subsection_Red1},
\ref{subsection_Red4}). The commutativity of the piece inside the
pentagon with the vertices $Z,Y,Y\quo L, X\quo G^\circ,X$ follows
directly from the definitions of the $L$-varieties
$Y,\underline{R}_Z,R$.

It remains to check that the piece with the vertices
$\underline{R}_Z\quo L,\l\quo L,X\quo G^\circ,$ $ \g\quo G$ is
commutative. Since the morphisms $Y\rightarrow \underline{R}_Z\quo
L$ and $X\rightarrow X\quo G^\circ$ are dominant, it is enough  to
prove the commutativity of the rectangle with vertices $Y,X,\l\quo
L,\g\quo G$. It is a direct consequence of
$\mu_{L,Y}=\mu_{G,X}|_{Y}$. \end{proof}

In the sequel we suppose that all morphisms between the varieties
from diagram (\ref{eq:5.7:1}) are the morphisms of this diagram.

\begin{lemma}\label{lemma:5.7.2}
The morphism $\widehat{\pi}\quo L':Y\quo L'\rightarrow R\quo L'$ is
birational.
\end{lemma}
\begin{proof}
The morphism $\widehat{\pi}\quo L'$ is dominant because so is
$\widehat{\pi}$.  Recall that $Y^{L'}\cong Y\quo L'$
(Corollary~\ref{corollary:4.2.3}). Under this identification,
$\widehat{\pi}\quo L'=\pi_{L',R}\circ \widehat{\pi}|_{Y^{L'}}$. Let
us note that $\pi_{P_0,X}|_{Y^{L'}}=\pi_{L',X\quo
P_u}\circ\nu\circ\widehat{\pi}|_{Y^{L'}}$. It follows from
Proposition~\ref{proposition:5.5.1} that there exists an open subset
$Y^1\subset Y^{L'}\cong Y\quo L'$ such that  for any $y_1\neq y_2\in
Y^1$ with $\widehat{\pi}(y_1)=\widehat{\pi}(y_2)$ there is $g\in
N_{G^\circ}(L,Y)\setminus L$ such that $y_1=gy_2$. But
$\widehat{\pi}(y_1)=\widehat{\pi}(y_2)$ implies
$\mu_{L,Y}(y_1)=\mu_{L,Y}(y_2)$. Since $\mu_{L,Y}(y_1)_s\subset
\z(\l)\cap\l^{pr}$, we see that $Z_{G^\circ}(\mu_{L,Y}(y_1))_s=L$.
Therefore $g\in L$. Thence $\widehat{\pi}(y_1)=\widehat{\pi}(y_2)$
yields $y_1=y_2$. \end{proof}

\begin{corollary}\label{corollary:5.7.3}
$\underline{\defe}_L(R)=\underline{\defe}_G(X),
\overline{\defe}_L(R)=\overline{\defe}_G(X)$.
\end{corollary}
\begin{proof}
In virtue of commutative diagram (\ref{eq:5.7:1}) and the fact that
$\widehat{\pi}$ is dominant,   $R$  is CN. Using
Corollary~\ref{corollary:4.2.3}, we have
$\underline{\defe}_L(R)=\underline{\defe}_{Z(L)^\circ}(R\quo L')$,
$\overline{\defe}_L(R)=\overline{\defe}_{Z(L)^\circ}(R\quo L')$.
Thanks to Lemma~\ref{lemma:5.7.2},
$\underline{\defe}_{Z(L)^\circ}(R\quo L' )=
\underline{\defe}_{Z(L)^\circ}(Y\quo L'),
\overline{\defe}_{Z(L)^\circ}(R\quo L'
)=\overline{\defe}_{Z(L)^\circ}(Y\quo{L'}).$ To complete the proof
apply Corollary \ref{corollary:4.2.3} to the action $L:Y$ and use
Proposition \ref{proposition:5.1.7}. \end{proof}

\begin{lemma}\label{lemma:5.7.4}
The subalgebra $A_Z^{L'}\subset A_Z$ is generated by $f|_Z,
H_\xi|_Z, f\in \K[X]^{P_0},\xi\in\z(\l)$.
\end{lemma}
\begin{proof}
Put $J=\Span_{A_Z}(H_\xi|_Z, \xi\in[\l,\l])$. Since
$\underline{R}_Z$ is  CN, $\mu_{L',\underline{R}_Z}^{-1}(0)\subset
\underline{R}_Z^{L'}$ and the restriction of
$\pi_{L',\underline{R}_Z}$ to $\mu_{L',\underline{R}_Z}^{-1}(0)$ is
a finite bijection (Corollary~\ref{corollary:4.2.3}). In particular,
the natural homomorphism $A_Z^{L'}\rightarrow A_Z/J$ is an
embedding. The image of this embedding coincides with
$(A_Z/J)^{L'}$.

For $f\in \K[X]^{P_u}, \xi\in \z(\l),$ denote by
$\overline{f},\overline{H}_\xi$ the image of $f,H_\xi$ in $A_Z/J$
under the natural epimorphism $A\rightarrow A_Z/J$. Clearly,
$\overline{f}, \overline{H}_\xi, f\in \K[X]^{P_u}, \xi\in \z(\l),$
generate $A_Z/J$. Note that $\overline{H}_{\xi}\in (A_Z/J)^{L'}$.
Let  $g\in (A_Z/J)^{L'}$. There exist $g_i\in \K[X]^{P_u},
\xi_i\in\z(\l)$ such that $g=\sum_i
\overline{g_i}\overline{H}_{\xi_i}$.  There is the natural
epimorphism of $L'$-modules $\K[X]^{P_u}\rightarrow \K[X]^{P_0}$.
Denote by $g_i^0$ the image of $g_i$ under this epimorphism. Since
$g,H_{\xi_i}$ are $L'$-invariant, $g=\sum_i
\overline{g_i^0}\overline{H}_{\xi_i}$. Hence the algebra
$(A_Z/J)^{L'}$ is generated by $\overline{H}_\xi,\xi\in\z(\l),
\overline{f}, f\in \K[X]^{P_0}$. It remains to recall that the
natural map $A_Z^{L'}\rightarrow (A_Z/J)^{L'}$ is an isomorphism.
\end{proof}

\begin{lemma}\label{lemma:5.7.5}
 The morphism $\nu\quo L':R\quo
L'\rightarrow X\quo P_0$ is finite.
\end{lemma}
\begin{proof}
The algebra $\K[R]^{L'}$  is integral over $A_Z^{L'}$. It remains to
check that  $A_Z^{L'}$ is integral over $\K[X]^{P_0}$. By
Lemma~\ref{lemma:5.7.4},
 it is enough to show that $H_{\xi}|_Z,\xi\in\z(\l),$ is integral
 over $\K[X]^{P_0}/(\K[X]^{P_0}\cap I_Z)$. This stems from diagram (\ref{eq:5.7:1}) because
 $H_{\xi}|_Z\subset \psi_{L,\underline{R}_Z}^*(\K[\l]^L)$ and $\K[\l]^L$
 is integral over $\K[\g]^G$.
 \end{proof}

Since the morphism $\widehat{\pi}:Y\rightarrow R$ is dominant, we
can identify $\K[R]$ with a subalgebra of $\K[Y]$.

\begin{corollary}\label{lemma:5.7.6}
The subalgebra $\K[R]^{L'}\subset \K[Y]^{L'}$ is the integral
closure of $\pi_{P_0,X}|_{Y}^*(\K[X]^{P_0})$.
\end{corollary}
\begin{proof}
$R\quo L'$ is a normal variety, the morphism $Y\quo L'\rightarrow
R\quo L'$ is birational (Lemma~\ref{lemma:5.7.2}), the morphism
$R\quo L'\rightarrow X\quo P_0$ is finite (Lemma~\ref{lemma:5.7.5}).
\end{proof}

\subsection{Proofs of Theorems~\ref{theorem:1.2.3},\ref{theorem:1.2.5}}\label{subsection_Red8}
We preserve the notation of
Subsections~\ref{subsection_Red5}-\ref{subsection_Red7}.

Recall that $Y$ is a Hamiltonian $N_0$-variety. Thus the morphism
$Y\quo L\rightarrow \l\quo L$ is $N_0/L$-equivariant. Further, the
morphism $Y\quo L\rightarrow X\quo G$ is the composition of
$\pi_{N_0/L,Y\quo L}$ and the open embedding $Y\quo
N_0\hookrightarrow X\quo G$.

\begin{lemma}\label{lemma:5.8.1}
The subalgebra $\K[R]^{L}\subset \K[Y]^L$ is  $N_0/L$-stable. The
morphism $\psi_{L,R}\quo L$ is $N_0/L$-equivariant, and the morphism
$R\quo L\rightarrow X\quo G$ is the quotient for the action
$N_0/L:R\quo L$.
\end{lemma}
\begin{proof}
The subalgebras $\K[X]^{G},\K[X]^{G^\circ}$ are embedded into
$\K[Y]^{L}$ via the restriction of functions to $Y$. By
Corollary~\ref{lemma:5.7.6}, $\K[R]^{L'}$ is the integral closure of
$\pi_{P_0,X}|_{Y}^*(\K[X]^{P_0})$. Therefore $\K[R]^L$ is the
integral closure of $\K[X]^{G^\circ}\subset \K[Y]^{L}$. Since
$\K[X]^{G^\circ}$ is integral over $\K[X]^G$, we obtain that
$\K[R]^L$ is  the integral closure of $\K[X]^G$ in $\K[Y]^L$. The
subalgebra $\K[R]^L\subset K[Y]^L$  is  $N_0/L$-stable because
$\K[X]^G\subset \K[Y]^{N_0}=(\K[Y]^L)^{N_0/L}$. Since
$\Quot(\K[Y]^{N_0})=\Quot(\K[X]^G)$, we have
$\Quot((\K[R]^L)^{N_0/L})=\Quot(\K[X]^G)$. Taking into account that
$\K[X]^G$ is integrally closed in $\Quot(\K[X]^G)$, we get
$(\K[R]^L)^{N_0/L}=\K[X]^G$. Let us recall that
$\widehat{\pi}:Y\rightarrow R$ is dominant. Since the morphisms
$\widehat{\pi}\quo L$ and $\psi_{L,Y}\quo L=\psi_{L,R}\quo
L\circ\widehat{\pi}\quo L$ are $N_0/L$-equivariant, so is
$\psi_{L,R}\quo L$. \end{proof}

We recall that $\a_{G,X}^{(Y)}=\overline{\im \mu_{Z(L)^\circ,Y\quo
L'}}$. Since $\widehat{\pi}:Y\rightarrow R$ is a dominant morphism
commuting with the moment maps, we have $\a_{G,X}^{(Y)}=\a_{L,R}$.

\begin{proposition}\label{lemma:5.8.2}
The subalgebra $\K[C_{L,R}]\subset \K[R]^L$ is
$W_{G,X}^{(Y)}$-stable and
$\K[C_{G,X}]=\K[C_{L,R}]^{W_{G,X}^{(Y)}}$. The following diagram is
commutative.

\begin{equation}\label{eq:5.8:1}
\begin{picture}(160,80)
\put(4,64){$R\quo L$}\put(4,24){$X\quo G$}
\put(60,64){$C_{L,R}$}\put(60,24){$C_{G,X}$}
\put(26,67){\vector(1,0){32}} \put(28,27){\vector(1,0){30}}
\put(14,60){\vector(0,-1){28}}\put(68,60){\vector(0,-1){28}}
\put(122,64){$\a_{G,X}^{(Y)}$}\put(114,24){$\a_{G,X}^{(Y)}/W_{G,X}^{(Y)}$}
\put(82,67){\vector(1,0){38}}\put(82,27){\vector(1,0){30}}
\put(132,60){\vector(0,-1){26}}\put(60,2){$\g\quo G$}
\put(68,22){\vector(0,-1){10}} \put(24,22){\vector(3,-1){34}}
\put(116,22){\vector(-3,-1){34}} 
\put(90,70){\scriptsize$\tau^2_{R,L}$}
\put(90,30){\scriptsize$\tau^2_{G,X}$}
\end{picture}
\end{equation}
\end{proposition}
\begin{proof}
We recall that $\K[X]^G=(\K[R]^L)^{W_{G,X}^{(Y)}}$
(Lemma~\ref{lemma:5.8.1}). It follows from commutative diagram
(\ref{eq:5.7:1}) that the subalgebra $\psi_{L,R}^*(\K[\l]^L)\subset
\K[X]^G$ is integral over $\psi_{G,X}^*(\K[\g]^G)$. Therefore
$\K[C_{L,R}]$ is the integral closure of $\K[C_{G,X}]$ in $\K[R]^L$.
In particular, $\K[C_{L,R}]$ is $W_{G,X}^{(Y)}$-stable. Since
$\K[C_{G,X}]$ is integrally closed in
$\K[X]^G=(\K[R]^L)^{W_{G,X}^{(Y)}}$, the equality
$\K[C_{G,X}]=\K[C_{L,R}]^{W_{G,X}^{(Y)}}$ holds.

Now we shall prove that the diagram (\ref{eq:5.8:1}) is commutative.
The only non-trivial thing here is the equality
$\tau^2_{G,X}=\tau^2_{L,R}/W_{G,X}^{(Y)}$. The latter is equivalent
to $\widehat{\psi}_{G,X}\quo G=(\widehat{\psi}_{L,R}\quo
L)/W_{G,X}^{(Y)}$. By the definition of $\widehat{\psi}_{G,X}$ (see
Subsection~\ref{subsection_Red2}),
$\widehat{\psi}_{G,X}|_{Y}=\widehat{\psi}_{N_0,Y}$. Therefore
$\widehat{\psi}_{G,X}\quo G|_{Y\quo N_0}=\widehat{\psi}_{N_0,Y}\quo
N_0$ (we recall that $Y\quo N_0$ is identified with an open subset
of $X\quo G$). Equivalently, $\widehat{\psi}_{G,X}\quo G|_{Y\quo
N_0}=(\widehat{\psi}_{L,Y}\quo L)/W_{G,X}^{(Y)}$. It remains to
recall that $\widehat{\psi}_{L,Y}\quo L=\widehat{\psi}_{L,R}\quo
L\circ \widehat{\pi}\quo L$ and that $\widehat{\pi}\quo L$ is a
$W_{G,X}^{(Y)}$-equivariant morphism. \end{proof}

\begin{proof}[of Theorem~\ref{theorem:1.2.3}]
Replacing $X$ with its normalization, we may assume that $X$ is normal.
The codimension of any irreducible component of a fiber of
$\psi_{L,R}$ in $R$ is not less than
$\underline{\defe}_L(R)=\underline{\defe}_G(X)$. Since the quotient
morphism $\pi_{L,R}$ is surjective, the same is true for any
irreducible component of a fiber of $\psi_{L,R}\quo L$. To complete
the proof it remains to apply Proposition~\ref{lemma:5.8.2}.
\end{proof}

\begin{proof}[of Theorem~\ref{theorem:1.2.5}]
The morphism $\widetilde{\psi}_{G,X}\quo G$ is equidimensional by
Theorem \ref{theorem:1.2.3}. Since $C_{G,X}$ is normal, it follows
that $\widetilde{\psi}_{G,X}\quo G$ is an open morphism
(see~\cite{Chevalley}). The equality
$\im\widetilde{\psi}_{G,X}=\im\widetilde{\psi}_{G,X}\quo G$ holds
because $\pi_{G,X}$ is surjective.

Put $Z=C_{L,R}, \tau=\tau^2_{L,R}:C_{L,R}\rightarrow \a_{L,R}\cong
\a_{G,X}^{(Y)}$. By Proposition~\ref{lemma:5.8.2}, $C_{G,X}\cong
Z/W_{G,X}^{(Y)}$, $\tau$ is $W_{G,X}^{(Y)}$-equivariant and
$\tau^2_{G,X}=\tau/W_{G,X}^{(Y)}$. Applying
Proposition~\ref{proposition:4.4.2} to the action $L:R$, we see that
$\tau$ is \'{e}tale in all points of $\im\widetilde{\psi}_{L,R}$. It
follows from commutative diagram~(\ref{eq:5.8:1})  that
$\im\widetilde{\psi}_{L,R}=\pi_{W_{G,X}^{(Y)},Z}^{-1}(\im\widetilde{\psi}_{G,X})$.
 \end{proof}

\subsection{The proof of Theorem~\ref{theorem:1.2.7}}\label{subsection_Red10}
The following proposition (at least, its first part) seems to be
quite standard.
\begin{proposition}\label{proposition:4.5.12}
Suppose $X$ is a generically symplectic normal affine irreducible
Hamiltonian $G$-variety.
\begin{enumerate}
\item The image of the embedding $\widetilde{\psi}_{G,X}^*:\K(C_{G,X})\rightarrow
\K(X)$ coincides with the center $\z(\K(X)^G)$ of the Poisson field
$\K(X)^G$.
\item Under the identification $\K(C_{G,X})\cong \z(\K(X)^G)$, the
equality
 $\K[\im\widetilde{\psi}_{G,X}]=\K[X]\cap
\z(\K(X)^G)$ holds.
\end{enumerate}
\end{proposition}
\begin{proof}
We identify $\K(C_{G,X})$ with
$\widetilde{\psi}_{G,X}^*(\K(C_{G,X}))$.

It follows from the definition of $C_{G,X}$ that $\K(C_{G,X})\subset
\K(X)^G$ and that $\K(C_{G,X})$ contains the subalgebra
$\psi_{G,X}^*(\K[\g]^G)$ and is algebraic over this subalgebra.
Since $\{H_\xi,f\}=L_{\xi_*}f=0$ for all $\xi\in\g, f\in \K(X)^G$,
the inclusion $\psi_{G,X}^*(\K[\g]^G)\subset \z(\K(X)^G)$ holds. The
uniqueness property for a lifting of a derivation yields
$\K(C_{G,X})\subset \z(\K(X))^G$.

Similarly to the proof of Satz 7.6 from \cite{Knop1}, one sees that
$\z_{\K(X)}(\K(X)^G)$ is algebraic over the subalgebra generated by
$H_\xi,\xi\in\g$. So $\z(\K(X)^G)=(\z_{\K(X)}(\K(X)^G))^G$ is
algebraic over $\widetilde{\psi}_{G,X}^*(\K[\g]^G)$ and thus also
over $\K(C_{G,X})$. To prove  assertion 1 of the proposition it
remains to show that $\K(C_{G,X})$ is algebraically closed in
$\K(X)^G$. This follows easily from the fact that $C_{G,X}$ is
integrally closed in $\K[X]^G$.

Proceed to assertion 2. Clearly,
$\K[\im\widetilde{\psi}_{G,X}]\subset \K[X]$. To prove the inclusion
$\K(C_{G,X})\cap \K[X]\subset \K[\im\widetilde{\psi}_{G,X}]$ note
that  the  pole locus of $\widetilde{\psi}_{G,X}^*(f)\in \K(X)$
coincides with the inverse image of the pole locus of $f\in
\K(C_{G,X})$.  \end{proof}

Now we consider a special class of Hamiltonian actions. Let $X$ be a
generically symplectic Hamiltonian $G$-variety.  One can show,
compare with~\cite{Vinberg}, Ch.2, $\S3$, Proposition 5, that the
following conditions are equivalent:
\begin{itemize}
\item[(a)] $\dim X=m_G(X)+\defe_G(X)$.
\item[(b)] The field $\K(X)^G$ is commutative with respect  to the Poisson bracket.
\end{itemize}

$X$ is called {\it coisotropic}, if it satisfies the equivalent
conditions (a),(b). The following statement follows immediately from
Proposition~\ref{proposition:4.5.12}.
\begin{corollary}\label{corollary:9.3}
Preserve the conventions of the previous proposition. If $X$ is
coisotropic, then $X\rightarrow \im\widetilde{\psi}_{G,X}$ is the
quotient morphism.
\end{corollary}

\begin{proof}[of Theorem~\ref{theorem:1.2.7}]
Let us prove that $\a_{G,X}^{(Y)}\subset \z(\l)$ is a linear
subspace. It follows from Lemma~\ref{lemma:3.3.7} that $0\in
\im\psi_{G,X}$. Since
$\psi_{G,X}=\tau^1_{G,X}\circ\widehat{\psi}_{G,X}$, we have
$0\in\im\tau^1_{G,X}$. It remains to apply
Proposition~\ref{proposition:4.4.1}.

Let us prove  that $\widetilde{\psi}_{G,X}$ is surjective. Let
$\lambda_0$ be such as in Lemma~\ref{lemma:3.3.7}. By the same
lemma, there is an action $\K^\times:C_{G,X}$ such that
$\widetilde{\psi}_{G,X}$ is $\K^\times$-equivariant and
$\lim_{t\rightarrow 0}t\lambda=\lambda_0$ for all $\lambda\in
C_{G,X}$. The image of $\widetilde{\psi}_{G,X}$ is
$\K^\times$-stable and contains $\lambda_0$. It follows from
Theorem~\ref{theorem:1.2.5} that $\im \widetilde{\psi}_{G,X}\subset
C_{G,X}$ is an open subset. We deduce that $\widetilde{\psi}_{G,X}$
is surjective. When $X$ is generically symplectic, we apply
Proposition~\ref{proposition:4.5.12} and obtain the equality
$\K[C_{G,X}]=\K[X]^G\cap \z(\K(X)^G)$.

It remains to prove that $\tau^2_{G,X}:C_{G,X}\rightarrow
\a_{G,X}/W_{G,X}$ is an isomorphism. The morphism $R\quo
L\rightarrow C_{G,X}$ from diagram~(\ref{eq:5.8:1}) is surjective
because  $\widetilde{\psi}_{G,X}\quo G$ is so. Taking into account
that the morphism $\widetilde{\psi}_{L,R}\quo L:R\quo L\rightarrow
C_{L,R}$ is $W_{G,X}^{(Y)}$-equivariant, we deduce from diagram
(\ref{eq:5.8:1}) that this morphism is also surjective.
Theorem~\ref{theorem:1.2.5} implies that
$\tau^2_{L,R}:C_{L,R}\rightarrow \a_{G,X}^{(Y)}$ is \'{e}tale. But
$\tau^2_{L,R}$ is finite and thus is an isomorphism. To complete the
proof it remains to apply Proposition~\ref{lemma:5.8.2}.
\end{proof}

\subsection{An example when $W_{G,X}^{(Y)}$ is not generated by reflections}\label{subsection_Red11}
In this subsection we construct an example of a conical symplectic
irreducible affine Hamiltonian $G$-variety $X$ with connected $G$
such that the image of $W_{G,X}^{(Y)}$ in $\GL(\a_{G,X}^{(Y)})$ is
not generated by reflections. Since $G$ is connected, the
homomorphism $W_{G,X}^{(Y)}\rightarrow \GL(\a_{G,X}^{(Y)})$ is an
embedding and we identify the group $W_{G,X}^{(Y)}$ with its image.

Put $G_0=\SL(2)\times \SL(2)$. Let $V_1,V_2$ (resp., $V_1',V_2'$) be
the two-dimensional irreducible modules over the first (resp., the
second) factor. Denote by $\gamma_V$ the linear automorphism of
$V:=V_1\oplus V_2\oplus V_1'\oplus V_2'$ given by the equality
$\gamma_V(v_1,v_2,v_1',v_2')=(v_1,-v_2,v_1',-v_2'), v_i\in V_i,
v_i'\in V_i'$. Put $G=\K^\times\times G_0,
\widetilde{X}=T^*\K^\times\times V$. Clearly, $\widetilde{X}$ is a
symplectic affine Hamiltonian $G$-variety. This is a conical
variety: the action $\K^\times:\widetilde{X}\cong G*_{G_0}(\K\oplus
V)$ is given by $t[g,(x,v)]=[g,(t^2x,tv)], g\in G, t\in \K^\times,
x\in K\cong T^*_1(\K^\times),v\in V$. Furthermore, $\widetilde{X}$
is coisotropic. Indeed, $m_G(\widetilde{X})=\dim G, \dim
\widetilde{X}=\dim G+\rank G$.

Let $\gamma_T$ denote the involution of $T^*\K^\times$ induced by
the left translation by $-1\in \K^\times$ and $\gamma$ the
involution of $\widetilde{X}$ given by $\gamma(x,y)=(\gamma_Tx,
\gamma_Vy), x\in T^*\K^\times, y\in V$. The involution $\gamma$ is
Hamiltonian and $\K^\times$-equivariant. Since $\gamma$ has no fixed
points, the variety $X=\widetilde{X}/\Z_2$,  where the non-unit
element of $\Z_2$ acts as $\gamma$, is smooth. So $X$ is a
symplectic conical coisotropic Hamiltonian  $G$-variety.

 The
variety $\widetilde{X}\quo G$ is isomorphic to $\A^3$. The action
$\Z_2:\widetilde{X}\quo G\cong \A^3$ is isomorphic to the linear
action of $\Z_2$ by the multiplication by matrices
$diag(1,\varepsilon,\varepsilon), \varepsilon=\pm 1$. We deduce that
$X\quo G\cong (\widetilde{X}\quo G)/\Z_2$ is not smooth. Using
Corollary~\ref{corollary:9.3} and Theorem~\ref{theorem:1.2.7}, we
see that $W_{G,X}^{(Y)}$ is not generated by reflections.

\section*{Index of notation}
As usual, if an algebraic group is denoted by a capital Latin
letter, we denote its Lie algebra by a corresponding small German
letter.
 \setlongtables
\begin{longtable}{p{4.5cm} p{12cm}}
$\langle\cdot,\cdot\rangle$& the pairing of elements of two dual to
each other vector spaces.
\\ $\# S$& the cardinality of a set $S$.
\\ $\partial_v$& the partial derivative in direction of a tangent vector $v$.
\\ $A^\times$& the group of invertible elements of an algebra $A$.
\\ $\operatorname{H}^0(X,V)$& the space of global  sections
of a vector bundle  $V$ over a variety $X$.
\\
 $G^\circ$& the connected component of unit of an algebraic group
$G$.\\ $(G,G)$ (resp., $[\g,\g]$)& the derived subgroup (resp., subalgebra) of a group $G$
(resp., the a Lie algebra $\g$).\\ $G*_HX$& the homogeneous bundle over
a homogeneous space $G/H$ with a fiber $X$.\\ $G_u$ (resp., $\g_u$)&
the unipotent radical of an algebraic group $G$ (resp.,
of an algebraic Lie algebra $\g$).\\
 $G_x$& the stabilizer of a point $x$
under an action of $G$.\\ $\g_x$& the annihilator of vector $x$ in a
module over a Lie algebra $\g$.
\\ $[g,x]$& the class of a pair $(g,x), g\in G, x\in X$ in the homogeneous bundle $G*_HX$.\\
$f|_Y$& the restriction of a map $f$ to a subset $Y$.
\\ $\g^\alpha$& the root subspace in a reductive Lie algebra $\g$
corresponding to  $\alpha\in\Delta(\g)$.
\\ $\im f$& the image of a map $f$.\\
$L_\xi$& the Lie derivative in direction of a vector field $\xi$.
\\ $m_G(X)$& the maximal dimension of an orbit for the action
of an algebraic group $G$ on a variety $X$.
\\  $N_G(H)$ (resp., $N_G(\h)$)&
the normalizer of a subgroup $H\subset G$ (resp., of a subalgebra
$\h$ of  $\g$) in a group $G$.\\
$\Quot(R)$& the quotient field of a domain $R$. \\
  $S_A(D)$& the symmetric algebra of
an $A$-module $D$.\\ $\Span_A(S)$& the $A$-submodule spanned by a
subset $S$ of some $A$-module.
\\ $\Spec(A)$& the affine scheme
corresponding to an algebra $A$.
\\ $U^\perp$& the orthogonal complement to a subspace $U\subset V$,  where $V$ is a vector space equipped
with a nondegenerate symmetric or skew-symmetric bilinear form.
\\ $v(f)$& the skew-gradient of a rational function  $f$ on a Poisson variety.
\\   $\overline{X}$& the closure of a subset of a variety with
respect to the Zariski topology.\\  $X^G$& the set of $G$-fixed
points for the action $G:X$.\\
 $X\quo
G$& the categorical quotient for the action of a reductive group $G$
on an affine variety $X$.\\ $X^{max}$& the subset of a Poisson
variety consisting of all  points  $x\in X^{reg}$ satisfying
$\rank P_x=\max_{y\in X^{reg}}\rank P_y$, where  $P$ is the Poisson
bivector of $X$.
\\ $X^{reg}$& the subset of smooth points of a variety $X$. \\ $Z(G)$ (resp., $\z(\g)$)&
the center of an algebraic group $G$ (resp.,  of a Lie algebra
$\g$).
\\ $\z_\g(\h)$& the centralizer of a subalgebra $\h$ in a Lie
algebra $\g$.
\\ $\alpha^\vee$& the dual root to a root $\alpha$.
\\ $\Delta(\g)$& the root system of a reductive Lie algebra
$\g$.
\\ $\xi_*$& the velocity vector field corresponding to an element
$\xi$ of a Lie algebra.
\\ $\xi_s$& the semisimple part of an element $\xi$ of a
reductive Lie algebra
\\
 $\pi_{G,X}$& the quotient morphism $X\rightarrow X\quo G$ for the action $G:X$.
\\
 $\varphi^*$& the morphism of algebras of  functions induced by a morphism $\varphi$ of varieties.\\
 $\varphi\quo G$& the morphism  $X\quo G\rightarrow Y\quo G$
induced by a  $G$-equivariant morphism $\varphi:X\rightarrow Y$.
\end{longtable}

\end{document}